\documentclass[reqno]{amsart}
\usepackage{amsfonts}
\usepackage{geometry,amsmath,amssymb}

\theoremstyle{plain}
 \newtheorem{theorem}{Theorem}[section]
 \newtheorem{proposition}{Proposition}[section]

\theoremstyle{definition}
 
 \newtheorem{definition}{Definition}[section]
\theoremstyle{remark}
 \newtheorem{remark}{Remark}[section]
 \numberwithin{equation}{section}

\begin{document}
\title[Stochastic fractional ...]
{Stochastic fractional evolution equations of order $\displaystyle 1< \alpha <2$ with generalized operators}
\author[Milo\v{s} Japund\v{z}i\'{c} and Danijela Rajter-\'{C}iri\'{c}]{$^{1}$Milo\v{s} Japund\v{z}i\'{c} and $^{2}$Danijela Rajter-\'{C}iri\'{c}}
\address{$^{1}$Faculty of Mathematics and Computer Sciences, Alfa BK University, Bulevar maršala Tolbuhina 8, 11000 Belgrade, Serbia \\
{\small\it e-mail:\mbox{ }milos.japundzic@gmail.com}}
\address{$^{2}$ Department of Mathematics and Informatics, Faculty of Science, University of Novi Sad, Trg Dositeja
Obradovi\'{c}a 4, 21000 Novi Sad, Serbia \\
{\small\it e-mail:\mbox{ }rajter@dmi.uns.ac.rs}}
\date{}

\begin{abstract}
We consider the Cauchy problem for stochastic fractional evolution equations with Caputo time fractional derivative of order $1<\alpha<2$ and space variable coefficients on an unbounded domain. The space derivatives that appear in the equations are of integer or fractional order such as the left and the right Liouville fractional derivative as well as the Riesz fractional derivative. To solve the problem we use generalized uniformly continuous solution operators. We obtain the unique solution within a certain Colombeau generalized stochastic process space. In our solving procedure, instead of the originate problem we solve a certain approximate problem, where operators of the original and the approximate problem are $L^2$-associated. Finally, application of the theory in solving stochastic time and time-space fractional wave equation is shown.

\medskip


\smallskip

\end{abstract}

\maketitle

\section{Introduction}

Wave phenomena arise in numerous real-life systems, making wave-type equations widespread and studied by many authors. Wave equations with standard (non-fractional) derivatives are sometimes good enough models, but in cases where long-term temporal memory or anomalous spatial diffusion is present, they fail to capture the complexity of the system. In these situations, fractional wave equations provide a more accurate framework for modeling. The choice of derivatives in the equation that will be taken as fractional depends on the properties of the system under consideration. When modeling wave phenomena with only temporal memory or history-dependent behavior in the medium, time-fractional wave equations are used (the fractional derivative is in time only, usually the Caputo derivative). For example, this is the case when the media has viscoelastic properties, such as biological tissues or polymers. On the other hand, space fractional or time-space fractional wave equations are good models for wave phenomena where spatial nonlocality or spatial memory occurs. Fractional wave equation are used when considering non-Newtonian fluids in porous media. For example, in oil recovery or groundwater flows, the porosity and fracture networks lead to spatially variable viscosity, which could be modeled with nonlocal fractional derivatives. Also these equations are convenient to use when considering fluid-structure interactions in fractured materials. In this case, both the elastic behavior (wave propagation) and fluid flow (with varying viscosity) need to be modeled together. This is typical in fractured rock or geophysical models.

In many physical, biological, and engineering systems, uncertainty and random perturbations play a significant role. For example, external fluctuations of the environment or inherent stochasticity in the medium may occur. In such cases, stochastic fractional wave equations are used. In these equations a stochastic source term and/or stochastic initial data appear.  External random inputs acting throughout the evolution of the system are represented by a stochastic source term in the equation (often modeled as time-space white noise). On the other hand, uncertainty in the initial state of the system is modeled by stochastic initial data in the equation. For example, this is appropriate when the initial displacement or wave velocity is not known deterministically. Studying fractional stochastic wave equations requires a function space framework capable of handling both nonlocal differential operators and low-regularity stochastic data. This motivates the use of fractional Sobolev spaces which provide the natural setting. These spaces interpolate between classical Sobolev and $L^2$ spaces and are well-suited to both fractional Laplacians and random data.

 Consder stochastic time and time-space fractional wave equations driven by a certain stochastic process $P$ (stochastic source term), for instance white noise process, and with stochastic initial data:
\begin{eqnarray}
& & \label{fwe-nonfrac}
  ^{C}{\mathcal{D}}_{t}^{\alpha}U(x,t) = \lambda(x)\partial_{x}^{2}U(x,t)+f(U(x,t))+P(x,t), \,\, (x,t)\in \mathbb{R}\times \mathbb{R}_{+}, \\
 & &  U(x,0)= U_{0}(x), \;\;\; U_{t}(x,0)= 0, \nonumber
 \end{eqnarray}
and
 \begin{eqnarray}
& & \label{fwe}
  ^{C}{\mathcal{D}}_{t}^{\alpha}U(x,t) = \lambda(x){\mathcal{D}}_{x}^{\beta}U(x,t)+f(U(x,t))+P(x,t), \,\, (x,t)\in \mathbb{R}\times \mathbb{R}_{+}, \\
 & &  U(x,0)= U_{0}(x), \;\;\; U_{t}(x,0)= 0. \nonumber
\end{eqnarray}
In equations above $^{C}{\mathcal{D}}_{t}^{\alpha}$ is Caputo time fractional derivative of order $ 1<\alpha<2,$ while ${\mathcal{D}}_{x}^{\beta}$ is space fractional derivative of order $1<\beta\leq2,$ such as the left and the right Liouville fractional derivative as well as the Riesz fractional derivative. Term $P(x,t)$ is a certain space-time generalized stochastic process. Initial data $U_0$ is a certain stochastic process.

The problems (\ref{fwe-nonfrac}) and (\ref{fwe}) can be rewritten in more general form:
 \begin{eqnarray}\label{fraceq}
\nonumber
  ^{C}{\mathcal{D}}_{t}^{\alpha}U(t) &=& AU(t)+f(U)+ P(\cdot,t), \,\, t>0, \\
  U(0)&=& U_{0}, \;\;\; U_{t}(0)= 0,
\end{eqnarray}
where, in both cases, the operator $A$  admits integer or fractional order space derivatives as well as variable coefficients depending on $x.$ So, instead of solving wave equations listed above, we will solve more general problems and (\ref{fraceq}) and then apply it to (\ref{fwe-nonfrac}) and (\ref{fwe}).

We solve the problem above in the framework of the Colombeau theory of generalized functions since our intention is to treat these equations by using an operator approach in contrast to commonly used numerical methods, that is, applying the solution operators as a generalization of semigroup of operators. For the Colombeau theory in general we refer, for example, to \cite{Biagoni}, \cite{Colombeau}, \cite{Nedeljkovb} and
\cite{Oberguggenberger}.

In our solving procedure, instead of the equation from the originate problem (\ref{fraceq}) we consider an approximate equation
\begin{equation}\label{fraceqaprox}
  ^{C}\mathcal{D}_{t}^{\alpha}U(t)=\widetilde{A}U(t)+f(U)+P(\cdot,t),
\end{equation}
with the same initial conditions, where $\widetilde{A}$ is a generalized linear bounded operator associated (in certain sense) to the originate operator $A$. In case of equations wave equations (\ref{fwe-nonfrac}) and (\ref{fwe}) one can obtain the approximate operator $\tilde A$ for a given (integer or fractional) differential operator $A$ by regularizing the derivative. The regularization is necessary in order to transform unbounded differential operators into bounded operators. To solve (\ref{fraceqaprox}) we use the theory of generalized uniformly continuous solution operators generated by $\widetilde A$, introduced and developed in \cite{EJDE}.

\section{Fractional derivatives and fractional Sobolev spaces}
%In this section we recall definitions of fractional derivatives with respect to time variable and give some useful estimates for fractional %derivatives and Mittag-Leffler function that we will use later.

%Recall first, the Fourier
%transform of a function $u\in L^{2}(\mathbb{R}),$
%denoted by $\mathcal{F}u\equiv \widehat{u},$ is defined as
%$$(\mathcal{F}u)(\xi)=\int\limits_{-\infty}^{\infty}e^{-ix\xi}u(x)dx,$$ while the inverse Fourier
%transform is given by $$u(x)\equiv
%(\mathcal{F}^{-1}(\mathcal{F}u))(x)=\frac{1}{2\pi}\int\limits_{-\infty}^{\infty}e^{ix\xi}\widehat{u}(\xi)d\xi.$$

\subsection{Time fractional derivatives}

The Caputo fractional derivative of order $\alpha,$ $m-1<\alpha\leq m,$ $m\in\mathbb{N},$ has the form (see, for example, \cite{Atanackovic}, \cite{Atanackovic1}, \cite{Kilbas}, \cite{Podlubny}, \cite{Samko})
$$ ^{C}\mathcal{D}_{t}^{\alpha}f(t)=J_{t}^{m-\alpha}f^{(m)}(t),\; \; \; \mbox{ where } \;
  J_{t}^{\alpha}f(t)=\frac{1}{\Gamma(\alpha)}\int_{0}^{t}(t-\tau)^{\alpha-1}f(\tau)d\tau,\; \; \alpha\geq0
$$
with $J^{0}=I,$ $I$ is identity operator. The Riemann-Liouville fractional derivative of order $\alpha,$ $m-1<\alpha\leq m,$ $m\in\mathbb{N},$ is given by $$  ^{RL}\mathcal{D}_{t}^{\alpha}f(t)=\frac{d^{m}}{dt^{m}}J_{t}^{m-\alpha}f(t).$$ For an absolutely continuous function $f$ and $1<\alpha<2,$ one has $J_{t}^{\alpha \;C}\mathcal{D}_{t}^{\alpha}f(t)=f(t)-f(0)-f^{'}(0)t$.

The following holds for a Riemann-Liouville fractional derivative:
\begin{proposition}\label{28} (\cite{Umarov}, Lemma 2.1) For all $\alpha\in (m-1,m]$ and $\beta\geq0$ the relation holds:
\begin{equation*}
  J_{t}^{\beta+\alpha}f(t)=J_{t}^{\beta+m \;RL}\mathcal{D}_{t}^{m-\alpha}f(t).
\end{equation*}
\end{proposition}

%Then for an absolutely continuous function $f$, the Caputo derivative in the case of $0<\alpha<1,$ can be written as (see \cite{Podlubny})
%$$^{C}\mathcal{D}_{t}^{\alpha}f(t)=\frac{d}{dt}f(t)*g_{1-\alpha}(t),$$
%where $*$ is the convolution operator given by
%$$f(t)*g(t)=\int_{-\infty}^{\infty}f(\tau)g(t-\tau)d\tau.$$
%Also, the Fourier transform of Caputo fractional derivative is given by %$$\widehat{^{C}\mathcal{D}_{t}^{\alpha}u}(\xi)=(i\xi)^{\alpha}\widehat{u}(\xi).$$
%Therefore, similar to Proposition 4 in \cite{Japundzic} we have the approximate Leibnitz rule for Caputo %fractional derivative as stated in the next proposition:
%\begin{prop}\label{26}
%Let $\alpha\in (0,1)$ and suppose that for all $t\geq0:$ $f(\cdot,t)\in L^{\infty}(\mathbb{R}),$ %$^{C}\mathcal{D}_{t}^{\alpha}f(\cdot,t) \in
%L^{2}(\mathbb{R}),$ $g(\cdot,t) \in L^{\infty}(\mathbb{R}),$
%$^{C}\mathcal{D}_{t}^{\alpha}g(\cdot,t) \in L^{2}(\mathbb{R}).$ Then one gets %$^{C}\mathcal{D}_{t}^{\alpha}(f(\cdot,t)g(\cdot,t)) \in L^{2}(\mathbb{R})$ and
%\begin{eqnarray*}
  %&&  \|^{C}\mathcal{D}_{t}^{\alpha}(f(\cdot,t)g(\cdot,t))\|_{L^{2}} \\
  %&&  \;\;\;\; \leq C\|f(\cdot,t)\|_{L^{\infty}}\|^{C}\mathcal{D}_{t}^{\alpha}g(\cdot,t)\|_{L^{2}}
  %+C\|^{C}\mathcal{D}_{t}^{\alpha}f(\cdot,t)\|_{L^{2}}\|g(\cdot,t)\|_{L^{\infty}}.
%\end{eqnarray*}
%\end{prop}

\subsection{Mittag-Leffler function}
The two-parameter Mittag-Leffler function $E_{\alpha,\beta}$ is given by $$E_{\alpha,\beta}(z)=\sum_{n=0}^{\infty}\frac{z^{n}}{\Gamma(\beta+n\alpha)},\; z\in \mathbb{C}, \; \alpha>0,\; \beta\in \mathbb{C}.$$ When $\beta=1$ we shortly write $E_{\alpha,1}(z)\equiv E_{\alpha}(z).$

For more details about Mittag-Leffler function we refer to \cite{EJDE}. Here we only state proposition that we use in proving procedures.

\begin{proposition}
Let $0<\alpha<2$ and $\beta>0$. Then
\begin{equation}\label{19}
  E_{\alpha,\beta}(\omega t^{\alpha})\leq C_{\alpha,\beta}(1+\omega^{(1-\beta)/\alpha})(1+t^{1-\beta})exp(\omega^{1/\alpha}t), \;\; \omega\geq0, \;\; t\geq0.
\end{equation}
\end{proposition}

The Cauchy problem (\ref{fraceq}) with Caputo fractional derivatives and nonlinear part $f(x,t)$ was considered in \cite{Umarov}, but in some special spaces of $L^{p}$ functions whose Fourier transforms are compactly supported in a some domain $G,$ and the following result was obtained:

\begin{proposition}\label{29-U} (\cite{Umarov}, Fractional Duhamel principle)\\
The solution of the Cauchy problem (\ref{fraceq}) is a given by
\begin{equation*}
  u(t)=E_{\alpha}(t^{\alpha}A)u_{0} + \int_{0}^{t}JE_{\alpha}((t-\tau )^{\alpha}A)^{RL}\mathcal{D}_{\tau}^{2-\alpha}f(\tau)d\tau.
\end{equation*}
\end{proposition}

\subsection{Space fractional derivatives and fractional Sobolev spaces}

Let $u\in C_{0}^{\infty}(\mathbb{R})$ and $m-1<\alpha<m,$ where $m\in\mathbb{N}$. Then the left Liouville fractional derivative of
order $\alpha$ on the whole axis $\mathbb{R}$
(see \cite{Kilbas}) is given by
$$(D_{+}^{\alpha}u)(x)=\frac{1}{\Gamma(m-\alpha)}\left(\frac{d}{dx}\right)^{m}\int\limits_{-\infty}^{x}\frac{u(\xi)}{(x-\xi)^{\alpha-m+1}}d\xi.$$
Since
$C_{0}^{\infty}(\mathbb{R})$ is a dense subset of
$H^{\alpha}(\mathbb{R}),$ this definition can be extended to a continuous linear map from $H^{\alpha}(\mathbb{R})$
into $L^{2}(\mathbb{R}).$ This extension of the left Liouville fractional derivative we denote by $\mathcal{D}_{+}^{\alpha},$ and simply
 call the Liouville fractional derivative.

We use the following definition of the fractional Sobolev space of order $\alpha$, $0<\alpha<1$:
 $$H^{\alpha}(\mathbb{R})=\{u\in L^{2}(\mathbb{R}): \mathcal{D}_{+}^{\alpha}u \in L^{2}(\mathbb{R})\},$$
with the norm $\|\cdot\|$ given by $
  \|u\|_{H^{\alpha}}=(\|u\|_{L^{2}}^{2}+\|\mathcal{D}_{+}^{\alpha}u\|_{L^{2}}^{2})^{\frac{1}{2}}.
$

Similarly, one can define fractional Sobolev space of order $\alpha$, $1<\alpha <2$, as:
$$
H^{\alpha}(\mathbb{R})=\{ u\in L^{2}(\mathbb{R}): \partial_{x}u\in L^{2}(\mathbb{R}), \mathcal{D}_{+}^{\alpha}u\in L^{2}(\mathbb{R})\} ,
$$
with the norm $\|\cdot\|$ given by $\|u\|_{H^{\alpha}}=(\|u\|_{L^{2}}^{2}+\|\partial_{x}u\|_{L^{2}}^{2}+\|\mathcal{D}_{+}^{\alpha}u\|_{L^{2}}^{2})^{\frac{1}{2}}.$

Note that, for $u\in H^{\alpha}(\mathbb{R}),$ $\alpha \geq 0$, it follows (see Corollary 3.2.9. in \cite{Webb}) that the Fourier transform of the Liouville fractional derivative of order $\alpha$ is $
   \|\mathcal{D}_{+}^{\alpha}u\|_{L^{2}}\leq\|u\|_{H^{\alpha}}.
$

Proofs of the propositions listed below can be found in \cite{Japundzic}.

\begin{proposition} Let $f\in C^{1}(\mathbb{C})$ be a Lipschitz function on $\mathbb{R}$ such that $f(0)=0$ and
suppose that $u\in H^{\alpha}(\mathbb{R}),
\hspace{1mm}0<\alpha<1.$ Then $\mathcal{D}_{+}^{\alpha}f(u)\in L^{2}$ and \begin{equation*}
  \|\mathcal{D}_{+}^{\alpha}f(u)\|_{L^{2}} \leq C\|f^{'}(u)\|_{L^{\infty}}\|\mathcal{D}_{+}^{\alpha}u\|_{L^{2}}.
\end{equation*}
\end{proposition}

\begin{proposition}\label{23} Let $\alpha\in (0,1)$ and suppose that $f\in L^{\infty},$ $\mathcal{D}_{+}^{\alpha}f \in
L^{2},$ $g \in L^{\infty},$
$\mathcal{D}_{+}^{\alpha}g \in L^{2}.$ Then $\mathcal{D}_{+}^{\alpha}(fg) \in L^{2}$ and
$$\|\mathcal{D}_{+}^{\alpha}(fg)\|_{L^{2}} \leq
C\|f\|_{L^{\infty}}\|\mathcal{D}_{+}^{\alpha}g\|_{L^{2}}+C\|\mathcal{D}_{+}^{\alpha}f\|_{L^{2}}\|g\|_{L^{\infty}}.$$
\end{proposition}

\begin{proposition}\label{190}
Suppose that $u\in H^{\alpha}(\mathbb{R}),
\hspace{1mm}1<\alpha<2$ and let $f\in C^{1}(\mathbb{C})$ be a Lipschitz function on $\mathbb{R}$ with bounded fractional derivatives $D_{+}^{\alpha}f$ such that $f(0)=0.$
 Then $\mathcal{D}_{+}^{\alpha}f(u)\in L^{2}$ and \begin{equation*}
  \|\mathcal{D}_{+}^{\alpha}f(u)\|_{L^{2}} \leq C\|f'(u)\|_{L^{\infty}}\|\mathcal{D}_{+}^{\alpha}u\|_{L^{2}}+C\|\mathcal{D}_{+}^{\alpha-1}(f'(u))\|_{L^{\infty}}\|u\|_{H^{\alpha}}.
\end{equation*}
\end{proposition}

\begin{proposition}\label{99}
Let $\alpha\in (1,2)$ and suppose that $f,g,\partial_{x}f,\partial_{x}g\in L^{\infty}$ and $\mathcal{D}_{+}^{\alpha-1}f,\mathcal{D}_{+}^{\alpha-1}g,\mathcal{D}_{+}^{\alpha}f,\mathcal{D}_{+}^{\alpha}g \in L^{2}$. Then $\mathcal{D}_{+}^{\alpha}(fg) \in L^{2}$ and
\begin{eqnarray*}
\nonumber
  \|\mathcal{D}_{+}^{\alpha}(fg)\|_{L^{2}} &\leq& C\|\partial_{x}f\|_{L^{\infty}}\|\mathcal{D}_{+}^{\alpha-1}g\|_{L^{2}}
   +C\|\mathcal{D}_{+}^{\alpha}f\|_{L^{2}}\|g\|_{L^{\infty}} \\
   && + C\|f\|_{L^{\infty}}\|\mathcal{D}_{+}^{\alpha}g\|_{L^{2}}+ C\|\mathcal{D}_{+}^{\alpha-1}f\|_{L^{2}}\|\partial_{x}g\|_{L^{\infty}}.
\end{eqnarray*}
\end{proposition}

The similar conclusions are also valid for the right Liouville fractional derivative of order $\alpha$
$(m-1<\alpha<m,$ $m\in\mathbb{N})$, on the whole axis $\mathbb{R}$ (see
\cite{Kilbas}) given by
$$(D_{-}^{\alpha}u)(x)=\frac{1}{\Gamma(m-\alpha)}\left(-\frac{d}{dx}\right)^{m}\int\limits_{x}^{\infty}\frac{u(\xi)}{(\xi-x)^{\alpha-m+1}}d\xi.$$
Namely, in the same manner as for the left fractional derivative, it is possible to extend this definition to a
continuous linear map from
$H^{\alpha}(\mathbb{R})$ into $L^{2}(\mathbb{R}),$ while the Fourier transform of the right Liouville fractional
derivative of order $\alpha\geq0$ is
$$\widehat{\mathcal{D}_{-}^{\alpha}u}(\xi)=(-i\xi)^{\alpha}\hat{u}(\xi).$$

The Riesz fractional derivative of order $m-1<\alpha<m$, $m\in\mathbb{N}$, denoted by $\mathcal{D}_{x}^{\alpha},$ is defined by using the left
and right $\alpha th$ Liouville fractional  derivative as
$$\mathcal{D}_{x}^{\alpha}u(x)=-\frac{1}{2\cos \frac{\alpha\pi}{2}}(\mathcal{D}_{+}^{\alpha}u(x)+\mathcal{D}_{-}^{\alpha}u(x)).$$

\section{Uniformly continuous solution operators}\label{ucso}
In \cite{EJDE} the theory of uniformly continuous solution operators is introduced and developed. Here we briefly repeat the main definitions and notions.

Consider the Cauchy problem for the fractional evolution equation of order $\alpha$, $1<\alpha<2:$
\begin{equation}\label{1}
  ^{C}{\mathcal{D}}_{t}^{\alpha}u(t)=Au(t), \, t>0; \,u(0)=x,\,u_{t}(0)=0,
\end{equation}
where $^{C}{\mathcal{D}}_{t}^{\alpha}$ is the Caputo fractional derivative of order $\alpha,$ and $A$ is a closed linear operator densely defined in a Banach space $E$. The problem (\ref{1}) is well-posed iff the following Volterra integral equation
\begin{equation}\label{zez}
  u(t)=x+\int_{0}^{t}g_{\alpha}(t-\tau)Au(\tau)d\tau
\end{equation}
is well-possed, where $g_{\alpha}(t)$ is defined for $\alpha>0,$ by
\begin{equation*}
g_{\alpha}(t)=\left\{
  \begin{array}{ll}
    t^{\alpha-1}/\Gamma(\alpha), & \hbox{$t>0$,} \\
    0, & \hbox{$t\leq 0$.}
  \end{array}
\right.
\end{equation*}

\begin{definition}
A family $S_{\alpha}(t), \; t\geq0,$ of linear and bounded operators on Banach space $E$ is called a uniformly continuous solution operator for (\ref{1}) if the following conditions are satisfied:
\begin{description}
  \item[(i)] $S_{\alpha}(t)$ is a uniformly continuous function for $t\geq0$ and $S_{\alpha}(0)=I,$ where $I$ is identity operator on $E$.
  \item[(ii)] $AS_{\alpha}(t)x=S_{\alpha}(t)Ax,$ for all $\,x\in E,\, t\geq0.$
  \item[(iii)] $S_{\alpha}(t)x$ is a solution of (\ref{zez}) for all $x\in E,$ $\, t\geq0$.
\end{description}
\end{definition}

\begin{definition}\label{4}
The infinitesimal generator A of a uniformly continuous solution operator $S_{\alpha}(t),$ $\alpha>0,$ $t\geq0,$  for (\ref{1}) is defined by
$$
  Ax=\Gamma(1+\alpha)\lim_{t\downarrow 0}\frac{S_{\alpha}(t)x-x}{t^{\alpha}}, \mbox{for all } x\in E.
$$
\end{definition}

\begin{definition}(\cite{Bazhlekova})
The solution operator $S_{\alpha}(t)$ is called exponentially bounded if there exist constants $M\geq1$ and $\omega\geq0$  such that $$\|S_{\alpha}(t)\|\leq Me^{\omega t}, \; t\geq0.$$
\end{definition}

From Definition \ref{4} it follows that every solution operator has a unique infinitesimal generator. If $S_{\alpha}(t)$ is a uniformly continuous solution operator satisfying $\|S_{\alpha}(t)\|\leq Me^{\omega t},$ for some $M\geq1$ and $\omega\geq0,$ its infinitesimal generator is a bounded linear operator.

On the other hand, every bounded linear operator $A$ is the infinitesimal generator of a uniformly continuous solution operator given by $$S_{\alpha}(t)=E_{\alpha}(t^{\alpha}A)=\sum_{n=0}^{\infty}\frac{t^{n\alpha}A^{n}}{\Gamma(1+n\alpha)}, \;\alpha>0,\;t\geq 0. $$
For every $1<\alpha\leq2$ this solution operator is unique as asserted in the following theorem.

\begin{proposition}[\cite{EJDE}]
Let $S_{\alpha}(t),$ $1<\alpha\leq2$, $t\geq0,$ be a uniformly continuous solution operator satisfying $\|S_{\alpha}(t)\|\leq Me^{\omega t},$ for some $M\geq1$ and $\omega\geq0.$ Then
\begin{description}
            \item[(i)] There exists a unique bounded linear operator $A$ such that $$S_{\alpha}(t)=E_{\alpha}(t^{\alpha}A), \; t\geq0.$$
            \item[(ii)] The operator $A$ in (i) is the infinitesimal generator of solution operator $S_{\alpha}(t).$
            \item[(iii)] For every $t\geq0$ holds $$^{C}{\mathcal{D}}_{t}^{\alpha}S_{\alpha}(t)=AS_{\alpha}(t)=S_{\alpha}(t)A.$$
          \end{description}
\end{proposition}

\section{Colombeau spaces}
The generalized operator which approximates the originate differential operator and generalized uniformly continuous solution operators belong to certain Colombeau spaces. Also, the space within which we solve the problem (\ref{fraceqaprox}) is a certain Colombeau space, too.

Let $(E, \|\cdot\|)$ be a Banach space and $\mathcal{L}(E)$ the space of all linear
continuous mappings from $E$ into $E$.

\begin{definition}
\label{Colombeau}
Let $m-1<\alpha<m,$ $m\in \mathbb{N}$. $\mathcal{S}E_{M}^{\alpha,m}([0,\infty):\mathcal{L}(E))$ is the space of nets $$(S_{\alpha})_{\varepsilon}:[0,\infty)\rightarrow
\mathcal{L}(E), \hspace{2mm} \varepsilon\in (0,1),$$
 with the following properties:
 \begin{description}
   \item[(i)] $(S_{\alpha})_{\varepsilon}(t)\in C^{m-1}([0,\infty):\mathcal{L}(E))\cap C^{m}((0,\infty):\mathcal{L}(E)).$
   \item[(ii)] $\lim_{t\rightarrow0^{+}}\left\|\frac{\frac{d^{m}}{dt^{m}}(S_{\alpha})_{\varepsilon}(t)}{t^{\alpha-m}}\right\|_{\mathcal{L}(E)}=C<+\infty.$
   \item[(iii)] For every $T>0$ there exist $N\in \mathbb{N},$ $M>0$ and $\varepsilon_{0}\in (0,1)$ such that
\begin{equation*}
  \sup_{t\in [0,T)} \left \|^{C}\mathcal{D}_{t}^{\gamma}(S_{\alpha})_{\varepsilon}(t)\right \|_{\mathcal{L}(E)}\leq
  M\varepsilon^{-N},
  \hspace{3mm}\varepsilon<\varepsilon_{0},\hspace{3mm} \gamma\in \{0,\ldots,m-1, \alpha\},
\end{equation*}
\begin{equation*}
  \sup_{t\in (0,T)} \left \|\frac{d^{m}}{dt^{m}}(S_{\alpha})_{\varepsilon}(t)\right \|_{\mathcal{L}(E)}\leq
  M\varepsilon^{-N},
  \hspace{3mm}\varepsilon<\varepsilon_{0}.
\end{equation*}
 \end{description}
 \end{definition}

Similarly we define the following space
\begin{definition}
\label{Colombeau1}
Let $m-1<\alpha<m,$ $m\in \mathbb{N}$. $\mathcal{S}N_{\alpha,m}([0,\infty):\mathcal{L}(E))$ is the space of nets
$$(N_{\alpha})_{\varepsilon}:[0,\infty)\rightarrow \mathcal{L}(E), \hspace{2mm} \varepsilon\in (0,1)$$ with the following properties:
 \begin{description}
   \item[(i)] $(N_{\alpha})_{\varepsilon}(t)\in C^{m-1}([0,\infty):\mathcal{L}(E))\cap C^{m}((0,\infty):\mathcal{L}(E)).$
   \item[(ii)] $\lim_{t\rightarrow0^{+}}\left\|\frac{\frac{d^{m}}{dt^{m}}(N_{\alpha})_{\varepsilon}(t)}{t^{\alpha-m}}\right\|_{\mathcal{L}(E)}=C<+\infty.$
   \item[(iii)] For every $T>0$ and $a\in\mathbb{R}$ there exist $M>0$ and $\varepsilon_{0}\in (0,1)$ such that
\begin{equation*}
  \sup_{t\in [0,T)} \left \|^{C}\mathcal{D}_{t}^{\gamma}(N_{\alpha})_{\varepsilon}(t)\right \|_{\mathcal{L}(E)}\leq
  M\varepsilon^{a},
  \hspace{3mm}\varepsilon<\varepsilon_{0},\hspace{3mm} \gamma\in \{0,\ldots,m-1, \alpha\},
\end{equation*}
\begin{equation*}
  \sup_{t\in (0,T)} \left \|\frac{d^{m}}{dt^{m}}(N_{\alpha})_{\varepsilon}(t)\right \|_{\mathcal{L}(E)}\leq
  M\varepsilon^{a},
  \hspace{3mm}\varepsilon<\varepsilon_{0}.
\end{equation*}
 \end{description}
 \end{definition}

In this paper, Caputo's fractional derivative in the problem of consideration is of order $1<\alpha <2$ and therefore, from now on, we will consider only that case. Hence, further we investigate spaces $\mathcal{S}E_{M}^{\alpha,2}([0,\infty):\mathcal{L}(E))$ and $\mathcal{S}N_{\alpha,2}([0,\infty):\mathcal{L}(E))$, although all the assertions we give can be extended for all $m\in {\mathbb N}$.

\begin{proposition}
\label{algebra}
The space $\mathcal{S}E_{M}^{\alpha,2}([0,\infty):\mathcal{L}(E))$ is an algebra with respect to composition of operators and the space $\mathcal{S}N_{\alpha,2}([0,\infty):\mathcal{L}(E))$ is an ideal of $\mathcal{S}E_{M}^{\alpha,2}([0,\infty):\mathcal{L}(E)).$
\end{proposition}

Now we can define a Colombeau-type space as a factor algebra by \begin{equation*}
  \mathcal{S}G_{\alpha,2}([0,\infty):\mathcal{L}(E))=\frac{\mathcal{S}E_{M}^{\alpha,2}([0,\infty):\mathcal{L}(E))}{\mathcal{S}N_{\alpha,2}([0,\infty):\mathcal{L}(E))}.
\end{equation*}
For every $1<\alpha<2$ elements of $\mathcal{S}G_{\alpha,2}([0,\infty):\mathcal{L}(E))$ will be denoted by $S=[(S_{\alpha})_{\varepsilon}],$
where $(S_{\alpha})_{\varepsilon}$ is a representative of the
class.

Similarly, one can define the following spaces:\\ $\mathcal{S}E_{M}(E)$ is the space of nets of linear
continuous mappings $$A_{\varepsilon}: E\rightarrow
E,\hspace{3mm} \varepsilon\in (0,1),$$ with the property that there exists constants $N\in \mathbb{N},$ $M>0$ and
$\varepsilon_{0}\in (0,1)$ such that
$$\|A_{\varepsilon}\|_{\mathcal{L}(E)}\leq M\varepsilon^{-N}, \hspace{3mm}\varepsilon<\varepsilon_{0}.$$
$\mathcal{S}N(E)$ is the space of nets of linear
continuous mappings $A_{\varepsilon}: E\rightarrow E,\hspace{3mm} \varepsilon\in (0,1),$ with the property that
for every $a\in \mathbb{R},$ there exist $M>0$
and $\varepsilon_{0}\in (0,1)$ such that $$\|A_{\varepsilon}\|_{\mathcal{L}(E)}\leq M\varepsilon^{a},
\hspace{3mm}\varepsilon<\varepsilon_{0}.$$ The Colombeau
space of generalized linear operators on $E$ is defined by
$$\mathcal{S}G(E)=\frac{\mathcal{S}E_{M}(E)}{\mathcal{S}N(E)}.$$ Elements of $\mathcal{S}G(E)$ will be
denoted by $A=[A_{\varepsilon}],$ where $A_{\varepsilon}$ is a representative of the class.

Finally, we introduce the Colombeau space within which we will solve (\ref{fraceqaprox}).\\
Let $1<\alpha<2,$ and $\frac{1}{2}<\beta<1$. $\mathcal{E}_{M}^{\alpha}([0,\infty):H^{\beta}(\mathbb{R}))$ is the
space of nets $$G_{\varepsilon}:[0,\infty)\times\mathbb{R}\rightarrow \mathbb{C}, \hspace{2mm} \varepsilon \in (0,1),$$
 with the following properties:
 \begin{description}
   \item[(i)] $G_{\varepsilon}(\cdot,\cdot)\in C^{1}([0,\infty):H^{\beta}(\mathbb{R}))\cap C^{2}((0,\infty):H^{\beta}(\mathbb{R})).$
   \item[(ii)] $\lim_{t\rightarrow0^{+}}\left\|\frac{\frac{d^{2}}{dt^{2}}G_{\varepsilon}(t,\cdot)}{t^{\alpha-2}}\right\|_{H^{\beta}}=C<+\infty.$
   \item[(iii)] For every $T>0$ there exist $M>0, N\in \mathbb{N}$ and
$\varepsilon_{0}>0$ such that
\begin{equation*}
  \sup_{t\in [0,T)} \left \|^{C}\mathcal{D}_{t}^{\gamma}G_{\varepsilon}(t,\cdot)\right \|_{H^{\beta}}\leq
  M\varepsilon^{-N},
  \hspace{3mm}\varepsilon<\varepsilon_{0},\hspace{3mm} \gamma\in \{0,1, \alpha\},
\end{equation*}
\begin{equation*}
  \sup_{t\in (0,T)} \left \|\frac{d^{2}}{dt^{2}}G_{\varepsilon}(t,\cdot)\right \|_{H^{\beta}}\leq
  M\varepsilon^{-N},
  \hspace{3mm}\varepsilon<\varepsilon_{0}.
\end{equation*}
 \end{description}
It is an algebra with respect to multiplication.\\
Similarly, for $1<\alpha<2,$ and $\frac{1}{2}<\beta<1,$ $\mathcal{N}_{\alpha}([0,\infty):H^{\beta}(\mathbb{R}))$ is the
space of nets $G_{\varepsilon}\in
\mathcal{E}_{M}^{\alpha}([0,\infty):H^{\beta}(\mathbb{R}))$ with the following properties:
 \begin{description}
   \item[(i)] $G_{\varepsilon}(\cdot,\cdot)\in C^{1}([0,\infty):H^{\beta}(\mathbb{R}))\cap C^{2}((0,\infty):H^{\beta}(\mathbb{R})).$
   \item[(ii)] $\lim_{t\rightarrow0^{+}}\left\|\frac{\frac{d^{2}}{dt^{2}}G_{\varepsilon}(t,\cdot)}{t^{\alpha-2}}\right\|_{H^{\beta}}=C<+\infty.$
   \item[(iii)] For every $T>0$ and $a\in\mathbb{R}$ there exist $M>0$ and $\varepsilon_{0}>0$ such that
\begin{equation*}
  \sup_{t\in [0,T)} \left \|^{C}\mathcal{D}_{t}^{\gamma}G_{\varepsilon}(t,\cdot)\right \|_{H^{\beta}}\leq
  M\varepsilon^{a},
  \hspace{3mm}\varepsilon<\varepsilon_{0},\hspace{3mm} \gamma\in \{0,1, \alpha\},
\end{equation*}
\begin{equation*}
  \sup_{t\in (0,T)} \left \|\frac{d^{2}}{dt^{2}}G_{\varepsilon}(t,\cdot)\right \|_{H^{\beta}}\leq
  M\varepsilon^{a},
  \hspace{3mm}\varepsilon<\varepsilon_{0}.
\end{equation*}
 \end{description}
The space $\mathcal{N}_{\alpha}([0,\infty):H^{\beta}(\mathbb{R}))$ is an ideal of $\mathcal{E}_{M}^{\alpha}([0,\infty):H^{\beta}(\mathbb{R}))$.

The quotient space
$$\mathcal{G}_{\alpha}([0,\infty):H^{\beta}(\mathbb{R}))=\frac{\mathcal{E}_{M}^{\alpha}([0,\infty):H^{\beta}(\mathbb{R}))}
{\mathcal{N}_{\alpha}([0,\infty):H^{\beta}(\mathbb{R}))}$$ is the corresponding Colombeau generalized function space
related to the Sobolev space $H^{\beta}$.

In the similar way, by omitting variable $t$, one can define spaces $\mathcal{E}_{M}^{\alpha}(H^{\beta}(\mathbb{R})),$
$\mathcal{N}_{\alpha}(H^{\beta}(\mathbb{R})),$ and
$\mathcal{G}_{\alpha}(H^{\beta}(\mathbb{R})).$

\section{Generalized uniformly continuous solution operators}
First, recall that every linear and bounded operator on Banach space $E$ is a closed and densely defined operator in $E$. Therefore, all results from Section \ref{ucso} continue to be valid in the case of linear and bounded operators on Banach space.

Instead of the Cauchy problem (\ref{1}) with closed and densely defined operator $A$, let us now consider fractional Cauchy problem given by
\begin{equation}\label{1a}
  ^{C}{\mathcal{D}}_{t}^{\alpha}u(t)=\widetilde{A}u(t), \, t>0; \,u(0)=x, \,u_{t}(0)=0,
\end{equation}
where $1<\alpha<2$ and $\widetilde{A}$ is a generalized linear bounded operator.

\begin{definition}
Let $1<\alpha<2$. $S_{\alpha}\in \mathcal{SG}_{\alpha,2}([0,\infty):\mathcal{L}(E))$  is called a Colombeau uniformly continuous solution operator for (\ref{1a}) if it has a representative $(S_{\alpha})_{\varepsilon}$ which is a uniformly continuous solution operator for (\ref{1a}) and for every $\varepsilon$ small enough.
\end{definition}

\begin{proposition}
Let $1<\alpha<2$ and let $(S_{\alpha})_{1\varepsilon}$ and $(S_{\alpha})_{2\varepsilon}$ be representatives of a generalized uniformly continuous solution operator $S_{\alpha}\in \mathcal{SG}_{\alpha,2}([0,\infty):\mathcal{L}(E)),$ with infinitesimal generators $\widetilde{A}_{1\varepsilon}$ and $\widetilde{A}_{2\varepsilon},$ respectively, for $\varepsilon$ small enough. Then $$\widetilde{A}_{1\varepsilon}-\widetilde{A}_{2\varepsilon}\in \mathcal{SN}(E).$$
\end{proposition}

\begin{definition}
$\widetilde A\in \mathcal{SG}(E)$ is called the infinitesimal generator of a Colombeau uniformly continuous solution operator $S_{\alpha}\in \mathcal{SG}_{\alpha,2}([0,\infty):\mathcal{L}(E))$, $1<\alpha<2$, if $\widetilde A_{\varepsilon}$ is the infinitesimal generator of the representative $(S_{\alpha})_{\varepsilon},$ for every $\varepsilon$ small enough.
\end{definition}

\begin{proposition}\label{20}
Let $1<\alpha<2$. Let $\widetilde A$ be the infinitesimal generator of a Colombeau uniformly continuous solution operator $S_{\alpha},$ and $\widetilde B$ the infinitesimal generator of a Colombeau uniformly continuous solution operator $T_{\alpha}.$ If $\widetilde A=\widetilde B$, then $S_{\alpha}=T_{\alpha}$.
\end{proposition}

Before we give the conditions under which generalized bounded operator is a generator of some generalized uniformly continuous solution operator some additionally properties will be defined.
\begin{definition}
Let $h_{\varepsilon}$ be a positive net satisfying $h_{\varepsilon}\leq \varepsilon^{-1}$. It is said that $\widetilde A\in \mathcal{SG}(E)$ is of $h_{\varepsilon}$-type if it has a representative $\widetilde A_{\varepsilon}$ such that $$\|\widetilde A_{\varepsilon}\|_{\mathcal{L}(E)}=\mathcal{O}(h_{\varepsilon}), \; \varepsilon\rightarrow0.$$
$G\in \mathcal{G}_{\alpha}([0,\infty):H^{\beta}(\mathbb{R}))$ is said to be of $h_{\varepsilon}$-type if it has a representative $G_{\varepsilon}$ such that $$\|G_{\varepsilon}\|_{H^{\beta}}=\mathcal{O}(h_{\varepsilon}), \; \varepsilon\rightarrow0.$$
\end{definition}
The following proposition holds for generalized operators.
\begin{proposition}\label{9}
Let $1<\alpha<2.$ Every $\widetilde A\in \mathcal{SG}(E)$ of $h_{\varepsilon}$-type, where $h_{\varepsilon}\leq C(\log1/\varepsilon)^{\alpha},$ is the infinitesimal generator of some generalized uniformly continuous solution operator $S_{\alpha} \in \mathcal{SG}_{\alpha,2}([0,\infty):{\mathcal{L}(E)}).$
\end{proposition}

Note that a Colombeau uniformly continuous solution opeator always possess an infinitesimal generator and it is unique. That follows from the fact that its representative is a classical uniformly continuous solution operator for which there exists a unique infinitesimal generator.

\section{Colombeau stochastic processes}
Now, we introduce some generalized function spaces that we will later need in order to define Colombeau generalized stochastic processes and space of generalized solutions to the problem of consideration.

Let $m-1<\beta <m$, $m\in {\mathbb N}$.

We denote $H^{\beta ,2}({\mathbb R})=\left \{ f \mbox{ such that }\partial^{\alpha}f \in L^{2}({\mathbb R}), \alpha \in \{0,\dots ,m-1, \beta \}  \right \}$.

We define the following spaces:

${\mathcal E}_{M}([0,\infty ): H^{\beta ,2} ({\mathbb R}))$ is the space of all
$$G_{\varepsilon}:(0,\infty )\times {\mathbb R}\mapsto {\mathbb C},\; G_{\varepsilon}(t,\cdot)\in H^{\beta ,2} ({\mathbb R}),\; \; \mbox{for every }t\in [0,\infty ),
$$
with the property that for every $T>0$ there exist $C>0$, $N\in {\mathbb N}$ and $\varepsilon_0 \in (0,1)$ such that
$$
\sup_{t\in [0,T)} \| \partial_t^{\alpha}G_{\varepsilon}(t,\cdot)\|_{H^{\beta ,2}({\mathbb R})}\leq C \varepsilon ^{-N}, \alpha \in \{0,1\}, \; \; \varepsilon < \varepsilon_0.
$$
We say that $\sup_{t\in [0,T)} \| \partial_t^{\alpha}G_{\varepsilon}(t,\cdot)\|_{H^{\beta ,2}({\mathbb R})}$ is moderate or that it has a moderate bound.

${\mathcal N}([0,\infty ): H^{\beta ,2} ({\mathbb R}))$ is the space of all $G_{\varepsilon}\in {\mathcal E}_{M}([0,\infty ): H^{\beta ,2} ({\mathbb R}))$ with the property that for every $T>0$ and $a\in {\mathbb R}$ there exist $C>0$ and $\varepsilon_0 \in (0,1)$ such that
$$
\sup_{t\in [0,T)} \| \partial_t^{\alpha}G_{\varepsilon}(t,\cdot)\|_{H^{\beta ,2}({\mathbb R})}\leq C \varepsilon ^{a}, \alpha \in \{0,1\}, \; \; \varepsilon < \varepsilon_0.
$$
We say that $\sup_{t\in [0,T)} \| \partial_t^{\alpha}G_{\varepsilon}(t,\cdot)\|_{H^{\beta ,2}({\mathbb R})}$ is negligible or that it has ${\mathcal N}$- bound.

Spaces ${\mathcal E}_{M}([0,\infty ): H^{\beta ,2} ({\mathbb R}))$ and ${\mathcal N}([0,\infty ): H^{\beta ,2} ({\mathbb R}))$ are algebras and ${\mathcal N}([0,\infty ): H^{\beta ,2} ({\mathbb R}))$ is an ideal of ${\mathcal E}_{M}([0,\infty ): H^{\beta ,2} ({\mathbb R}))$, so we can define the factor algebra
$$
{\mathcal G}([0,\infty ): H^{\beta ,2} ({\mathbb R}))=\frac {{\mathcal E}_{M}([0,\infty ): H^{\beta ,2} ({\mathbb R}))}{{\mathcal N}([0,\infty ): H^{\beta ,2} ({\mathbb R}))}
$$
which is called the algebra of $H^{\beta ,2}$-Colombeau generalized functions.

By omitting the variable $t$, one can similarly define the spaces ${\mathcal E}_{M}(H^{\beta ,2}({\mathbb R}))$, ${\mathcal N}(H^{\beta ,2}({\mathbb R}))$ and ${\mathcal G}(H^{\beta ,2}({\mathbb R}))$.

%\begin{definition}\label{GLinfty process}
%A ${\mathcal G}_{H^{2,\infty}}$-Colombeau generalized stochastic process on a probability space $\left ( \Omega, \Sigma, \mu \right )$ is a %mapping $U: \Omega \mapsto {\mathcal G}([0,\infty ): H^{2,\infty} ({\mathbb R}))$ such that there exists a function $\tilde{U}:(0,1)\times %[0,\infty )\times {\mathbb R} \times \Omega \mapsto {\mathbb R}$ with the following properties:
%\begin{itemize}
%\item[1)] For fixed $\varepsilon \in (0,1),\; (t,x,\omega )\mapsto \tilde{U}(\varepsilon ,t,x, \omega )$ is jointly measurable in $[0,\infty %)\times {\mathbb R} \times \Omega$.
%\item[2)] The mapping $\varepsilon \mapsto \tilde{U}(\varepsilon ,t,x, \omega )$ is en element of ${\mathcal E}_{M}([0,\infty ): H^{2,\infty} %({\mathbb R}))$ almost surely in $\omega \in \Omega$, and it is a representative of $U(\omega )$.
%\end{itemize}
%The algebra of ${\mathcal G}_{H^{2,\infty}}$-Colombeau generalized stochastic processes on $\Omega$ will be denoted by ${\mathcal %G}^{\Omega}([0,\infty ): H^{2,\infty} ({\mathbb R}))$.
%\end{definition}

\begin{definition}\label{GHbeta2 process}
A ${\mathcal G}_{H^{\beta, 2}}$-Colombeau generalized stochastic processes on a probability space $\left ( \Omega, \Sigma, \mu \right )$ is a mapping $U: \Omega \mapsto {\mathcal G}([0,\infty ): H^{\beta, 2} ({\mathbb R}))$ such that there exists a function $\tilde{U}:(0,1)\times [0,\infty )\times {\mathbb R} \times \Omega \mapsto {\mathbb R}$ with the following properties:
\begin{itemize}
\item[1)] For fixed $\varepsilon \in (0,1),\; (t,x,\omega )\mapsto \tilde{U}(\varepsilon ,t,x, \omega )$ is jointly measurable in $[0,\infty )\times {\mathbb R} \times \Omega$.
\item[2)] The mapping $\varepsilon \mapsto \tilde{U}(\varepsilon ,t,x, \omega )$ is en element of ${\mathcal E}_{M}([0,\infty ): H^{\beta, 2} ({\mathbb R}))$ almost surely in $\omega \in \Omega$, and it is a representative of $U(\omega )$.
\end{itemize}
The algebra of ${\mathcal G}_{H^{\beta, 2}}$-Colombeau generalized stochastic processes on $\Omega$ will be denoted by ${\mathcal G}^{\Omega}([0,\infty ): H^{\beta, 2} ({\mathbb R}))$.
\end{definition}

Finally, note that by omitting variable $t$ one can define space ${\mathcal G}^{\Omega}(H^{\beta, 2} ({\mathbb R}))$ as a space of mappings $U: \Omega \mapsto {\mathcal G}(H^{\beta, 2} ({\mathbb R}))$ such that there exists a function $\tilde{U}:(0,1)\times {\mathbb R} \times \Omega \mapsto {\mathbb R}$ with the following properties:
\begin{itemize}
\item[1)] For fixed $\varepsilon \in (0,1),\; (x,\omega )\mapsto \tilde{U}(\varepsilon ,x, \omega )$ is jointly measurable in ${\mathbb R} \times \Omega$.
\item[2)] The mapping $\varepsilon \mapsto \tilde{U}(\varepsilon ,x, \omega )$ is en element of ${\mathcal E}_{M}(H^{\beta, 2} ({\mathbb R}))$ almost surely in $\omega \in \Omega$, and it is a representative of $U(\omega )$.
\end{itemize}

\vspace{1.5mm}

In the sequel the variable ${\varepsilon}$ will be written as a subindex. Thus, $U_{\varepsilon}$ will always denote a representative of $U$.

\section{The existence-uniqueness result}

We start with a case:
\begin{eqnarray}\label{fraceq-nulti}
\nonumber
  ^{C}{\mathcal{D}}_{t}^{\alpha}u(t) &=& Au(t)+F(\cdot,t,u), \,\, t>0, \\
  u(0)&=& u_{0}, \;\;\; u_{t}(0)= 0.
\end{eqnarray}

Motivated by Proposition \ref{29-U} it follows that fractional Duhamel principle in the case of solution operator has the form:

\begin{proposition}\label{frac2}
The solution of the Cauchy problem (\ref{fraceq-nulti}) is a given by:
 \begin{equation}\label{36}
  u(t)=S_{\alpha}(t)u_{0}
+\int_{0}^{t}\;_{\tau}J_{t}S_{\alpha}(t-\tau)^{RL}\mathcal{D}_{\tau}^{2-\alpha}F( \cdot , \tau , u(\tau))d\tau,
\end{equation}
where $S_{\alpha}(t)$ is a solution operator generated by $A$.
\end{proposition}

\begin{remark}
In that case when $F(x, 0, u_{0})=0$ one can replace Riemann-Liouville derivative with Caputo derivative in (\ref{36}).
\end{remark}

Integral representation stated in the next proposition will often be used in proving some auxiliary results as well as in proving our main result. (For the proofs we  refer to \cite{EJDE}).

\begin{proposition}\label{intrepre1}
Let $1<\alpha<2$ and let $S_{\alpha}(t)$ be a solution operator generated by $A$. Then
\begin{eqnarray}\label{35}
\nonumber
 && \int_{0}^{t}\;_{\tau}J_{t}S_{\alpha}(t-\tau)^{RL}\mathcal{D}_{\tau}^{2-\alpha}f(u(\tau))d\tau\\
  &&\hspace{15mm}=\int_{0}^{t}(t-\tau)^{\alpha-1}E_{\alpha,\alpha}((t-\tau)^{\alpha}A)f(u(\tau))d\tau.
  \end{eqnarray}
\end{proposition}

Later, in proving certain results in Colombeau framework, we will need second order derivative of an integral appearing in the representation (\ref{35}).
\begin{proposition}
Let $1<\alpha<2$ and let $S_{\alpha}(t)$ be a solution operator generated by $A$. Then
\begin{eqnarray}\label{35c}
\nonumber
 && \frac{d}{dt^{2}}\int_{0}^{t}\;_{\tau}J_{t}S_{\alpha}(t-\tau)^{RL}\mathcal{D}_{\tau}^{2-\alpha}f(u(\tau))d\tau\\\nonumber
  &&\hspace{10mm}=J_{t}^{\alpha-1}\frac{d}{dt}f(u(t)) + \frac{f(u(0))t^{\alpha-2}}{\Gamma(\alpha-1)}\\ &&\hspace{10mm}+A\int_{0}^{t}(t-\tau)^{2\alpha-3}E_{\alpha,2\alpha-2}((t-\tau)^{\alpha}A)f(u(\tau))d\tau.
  \end{eqnarray}
\end{proposition}

We now take $E=H^{\beta}(\mathbb{R}),$ $\frac{1}{2}<\beta<1$. Instead of the Cauchy problem (\ref{fraceq}) with closed and densely defined operator $A$ on $H^{\beta}(\mathbb{R})$, we will consider fractional Cauchy problem given by
$$
  ^{C}{\mathcal{D}}_{t}^{\alpha}U(t)=\widetilde{A}U(t)+f(
  U)+P(\cdot,t), \, t>0; \,\,U(0)=U_{0},\,U_{t}(0)=0,
$$
where $\widetilde{A}$ is a generalized linear bounded operator $L^{2}-$associated with $A$, i.e., for every $u\in H^{\beta}(\mathbb{R})$, the following holds $$\|(A-\widetilde{A}_{\varepsilon})u\|_{L^{2}}\rightarrow0,\; \; \varepsilon\rightarrow0.$$

\begin{theorem}\label{14}
Let $1<\alpha<2$, $\frac{1}{2}<\beta<1,$ and $f\in C^{2}(\mathbb R)$ such that $f$ and $f'$ are Lipschitz on $\mathbb R$ satisfying $f(0)=f'(0)=0.$ Let Colombeau generalized processes $Q$ and $P$ be such that $Q\in {\mathcal G}^{\Omega}(H^{\beta} ({\mathbb R}))$ and $P\in {\mathcal G}^{\Omega}([0,\infty ): H^{\beta} ({\mathbb R}))$.
Suppose that operator $\widetilde{A}\in \mathcal{SG}(H^{\beta}(\mathbb{R}))$ is of $h_{\varepsilon}$-type, with $h_{\varepsilon}= o\Big((\log(\log1/\varepsilon)^{\alpha-1})^{\alpha}\Big)$, such that $\|\widetilde{A}_{\varepsilon}u_{\varepsilon}\|_{L^{2}}\leq h_{\varepsilon}\|u_{\varepsilon}\|_{L^{2}},$ for $u_{\varepsilon}\in H^{\beta}(\mathbb{R}).$
%and such that closed and densely defined operator $A$ from (\ref{1}) and $\widetilde{A}$ are $L^{2}-$associated, i.e. for every $u\in H^{1}(\mathbb{R})$, the following holds $$\|(A-\widetilde{A}_{\varepsilon})u\|_{L^{2}}\rightarrow0,\; \; \varepsilon\rightarrow0.$$

 Then, for every $1<\alpha<2$ and $\frac{1}{2}<\beta<1,$ a unique generalized solution $U\in \mathcal{G}^{\Omega}_{\alpha}([0,\infty):H^{\beta}(\mathbb{R}))$ to the Cauchy problem
\begin{eqnarray}\label{6b}
\nonumber
  ^{C}{\mathcal{D}}_{t}^{\alpha}U(t) &=& \widetilde{A}U(t)+f(U)+P(\cdot,t), \,\, t>0, \\
  U(0)&=& Q, \;\;\; U_{t}(0)= 0,
\end{eqnarray}
almost surely exists and it is represented by
 \begin{equation}\label{7b}
  U_{\varepsilon}(t)=(S_{\alpha})_{\varepsilon}(t)Q_{\varepsilon}
+\int_{0}^{t}\,_{\tau}J_{t}(S_{\alpha})_{\varepsilon}(t-\tau)^{RL}\mathcal{D}_{\tau}^{2-\alpha}(f(U_{\varepsilon}(\tau))+P_{\varepsilon}(\cdot,\tau))d\tau,
\end{equation}
where  $S_{\alpha}\in
\mathcal{SG}_{\alpha,2}([0,\infty):\mathcal{L}(H^{\beta}(\mathbb{R})))$  is a Colombeau uniformly continuous solution operator generated by $\widetilde{A}$.
\end{theorem}

\begin{remark} Note that a term $f(U_{\varepsilon}(\tau))+P_{\varepsilon}(\cdot,\tau)$ in (\ref{7b}) is deterministic (as a representative of a Colombeau stochastic process) and therefore it is a special case of $F$ in (\ref{fraceq-nulti}).
\end{remark}

\begin{proof}
Fix $1<\alpha<2$ and $\frac{1}{2}<\beta<1.$ Since from the assumption follows that the operator $\widetilde{A}$ is of $h_{\varepsilon}$-type, with $h_{\varepsilon}= o((\log\log1/\varepsilon)^{\alpha}),$ it is obvious that the operator $\widetilde{A}$ is the infinitesimal generator of a Colombeau solution operator $S_{\alpha}\in
\mathcal{SG}_{\alpha,2}([0,\infty):\mathcal{L}(H^{\beta}(\mathbb{R})))$ given by $S_{\alpha}(t)=E_{\alpha}(t^{\alpha}\widetilde{A})$ (see Proposition \ref{9}). Also, from (\ref{36}) we know that (\ref{7b}) represents a solution to (\ref{6b}).

 Let us show that this solution is an element of $\mathcal{G}^{\Omega}_{\alpha}([0,\infty):H^{\beta}(\mathbb{R})).$

 Fix $\omega \in \Omega$. The mapping $(x,t,\omega)\mapsto U_{\varepsilon}(x,t,\omega)$ is a jointly measurable in $(x,t)$ and $\omega$ for every fixed $\varepsilon$ which follows by applying the successive approximation method and the fact that a continuous mapping of a measurable function is also measurable.

 First, we show that the solution satisfies
\begin{equation}\label{26a}
  \lim_{t\rightarrow0^{+}}\left\|\frac{\frac{d}{dt^{2}}U_{\varepsilon}(t,\cdot)}{t^{\alpha-2}}\right\|_{H^{\beta}}=C<+\infty.
\end{equation}
Indeed, after differentation of (\ref{7b}) twice with respect to $t,$ using integral  representation (\ref{35c}) one gets
 \begin{eqnarray}\label{33c}
 \nonumber
   \frac{d}{dt^{2}}U_{\varepsilon}(t,\cdot) &=& \frac{d}{dt^{2}}(S_{\alpha}^{})_{\varepsilon}(t)Q_{\varepsilon}+J_{t}^{\alpha-1}\frac{d}{dt}(f(U_{\varepsilon}(t))+{P_\varepsilon}(\cdot,t)) + \frac{f(Q_{\varepsilon})+{P_\varepsilon}(\cdot,0)}{\Gamma(\alpha-1)}t^{\alpha-2} \\\nonumber
    && +\int_{0}^{t}(t-\tau)^{2\alpha-3}\widetilde{A}_{\varepsilon}E_{\alpha,2\alpha-2}((t-\tau)^{\alpha}\widetilde{A}_{\varepsilon})(f(U_{\varepsilon}(\tau))+{P_\varepsilon}(\cdot,\tau))d\tau.\\
 \end{eqnarray}
 Further we have
 \begin{eqnarray}\label{34b}
 \nonumber
   &&\|\frac{d}{dt^{2}}U_{\varepsilon}(t,\cdot)\|_{L^{2}}\\
    \nonumber  &&\leq\|\frac{d}{dt^{2}}(S_{\alpha}^{})_{\varepsilon}(t)Q_{\varepsilon}\|_{L^{2}}+\frac{1}{\Gamma(\alpha-1)}\int_{0}^{t}(t-\tau)^{\alpha-2}
  \left\|\partial_{\tau}\Big(f({U_{\varepsilon}(\tau))}+{P_\varepsilon}(\cdot,\tau)\Big)\right\|_{L^{2}}d\tau \\
   \nonumber
  &&\hspace{5mm} +\frac{\|f(Q_{\varepsilon})\|_{L^{2}}+\|{P_\varepsilon}(\cdot,0))\|_{L^{2}}}{\Gamma(\alpha-1)}t^{\alpha-2}\\\nonumber
  &&\hspace{5mm}+\int_{0}^{t}(t-\tau)^{2\alpha-3}\|\widetilde{A}_{\varepsilon}\|\,
  E_{\alpha,2\alpha-2}((t-\tau)^{\alpha}\|\widetilde{A}_{\varepsilon}\|)\,(\|f(U_{\varepsilon}(\tau))\|_{L^{2}}+\|{P_\varepsilon}(\cdot,\tau)\|_{L^{2}}) d\tau\\\nonumber
  &&\leq \|\frac{d}{dt^{2}}(S_{\alpha})_{\varepsilon}(t)Q_{\varepsilon}\|_{L^{2}}
  +\frac{1}{\Gamma(\alpha-1)}\int_{0}^{t}(t-\tau)^{\alpha-2}
  \|f'(U_{\varepsilon}(\tau))\|_{L^{\infty}}\|\partial_{\tau}U_{\varepsilon}(\tau)\|_{L^{2}}d\tau\\\nonumber
  &&\hspace{5mm}+\frac{1}{\Gamma(\alpha-1)}\int_{0}^{t}(t-\tau)^{\alpha-2}
  \|\partial_{\tau}{P_\varepsilon}(\cdot,\tau)\|_{L^{2}}d\tau+\frac{C\|Q_{\varepsilon}\|_{L^{2}}+\|{P_\varepsilon}(\cdot,0))\|_{L^{2}}}{\Gamma(\alpha-1)}t^{\alpha-2}\\
  &&\hspace{5mm}+\int_{0}^{t}(t-\tau)^{2\alpha-3}\|\widetilde{A}_{\varepsilon}\|\,
  E_{\alpha,2\alpha-2}((t-\tau)^{\alpha}\|\widetilde{A}_{\varepsilon}\|)
  \,\Big(C\|U_{\varepsilon}(\tau)\|_{L^{2}}+ \|{P_\varepsilon}(\cdot,\tau)\|_{L^{2}}\Big)d\tau.
 \end{eqnarray}
Knowing that for function $E_{\alpha,2\alpha-2}$ following estimate holds
 \begin{equation}\label{2alpha}
   E_{\alpha,2\alpha-2}(t^{\alpha}\|\widetilde{A}_{\varepsilon}\|)\leq C_{\alpha}(1+\|\widetilde{A}_{\varepsilon}\|^{\frac{3-2\alpha}{\alpha}})(1+t^{3-2\alpha})
   exp\Big(\|\widetilde{A}_{\varepsilon}\|^{\frac{1}{\alpha}}t\Big),
 \end{equation}
 setting
 \begin{equation*}
   \widetilde{M}_{t}=\sup_{\tau \in [0,t)}C_{\alpha}(1+\|\widetilde{A}_{\varepsilon}\|^{\frac{3-2\alpha}{\alpha}})
   exp\Big(\|\widetilde{A}_{\varepsilon}\|^{\frac{1}{\alpha}}\tau\Big),
 \end{equation*}
 the last integral in (\ref{34b}) can be estimated by
 \begin{eqnarray*}
    &&  \int_{0}^{t}(t-\tau)^{2\alpha-3}\|\widetilde{A}_{\varepsilon}\|\,
  E_{\alpha,2\alpha-2}((t-\tau)^{\alpha}\|\widetilde{A}_{\varepsilon}\|)
  \,\Big(C\|U_{\varepsilon}(\tau)\|_{L^{2}}+ \|{P_\varepsilon}(\cdot,\tau)\|_{L^{2}}\Big)d\tau. \\
    && \hspace{10mm}\leq \widetilde{M}_{t}(t+\frac{t^{2\alpha-2}}{2\alpha-2})\|\widetilde{A}_{\varepsilon}\|\sup_{\tau\in [0,t)}\Big(C\|U_{\varepsilon}(\tau)\|_{L^{2}}+ \|{P_\varepsilon}(\cdot,\tau)\|_{L^{2}}\Big),
 \end{eqnarray*}
 which implies
 \begin{equation*}
  \lim_{t\rightarrow0^{+}}\left\|\frac{\frac{d}{dt^{2}}U_{\varepsilon}(t,\cdot)}{t^{\alpha-2}}\right\|_{L^{2}}=C<+\infty.
\end{equation*}

  Also, after differentiation of (\ref{33c}) with respect to left Liouville fractional derivative $\mathcal{D}_{+}^{\beta},$ one gets
 \begin{eqnarray}\label{35b}
  \nonumber
   &&\|\mathcal{D}_{+}^{\beta}\frac{d}{dt^{2}}U_{\varepsilon}(t,\cdot)\|_{L^{2}}\leq\|\frac{d}{dt^{2}}(S_{\alpha}^{})_{\varepsilon}(t)\mathcal{D}_{+}^{\beta}Q_{\varepsilon}\|_{L^{2}}\\
    \nonumber
    &&\hspace{5mm}+\frac{1}{\Gamma(\alpha-1)}\int_{0}^{t}(t-\tau)^{\alpha-2}
  \Big(\left\|\mathcal{D}_{+}^{\beta}\partial_{\tau} f({U_{\varepsilon}(\tau)})\right\|_{L^{2}}+\left\|\mathcal{D}_{+}^{\beta}\partial_{\tau}{P_\varepsilon}(\cdot,\tau)\right\|_{L^{2}}\Big)d\tau \\
   \nonumber
   &&\hspace{5mm} +\frac{\|\mathcal{D}_{+}^{\beta}f(Q_{\varepsilon})\|_{L^{2}}+\|\mathcal{D}_{+}^{\beta}{P_\varepsilon}(\cdot,0))\|_{L^{2}}}{\Gamma(\alpha-1)}t^{\alpha-2}\\\nonumber
  &&\hspace{5mm}+\int_{0}^{t}(t-\tau)^{2\alpha-3}\|\widetilde{A}_{\varepsilon}\|\,
  E_{\alpha,2\alpha-2}((t-\tau)^{\alpha}\|\widetilde{A}_{\varepsilon}\|)\Big(\left\|\mathcal{D}_{+}^{\beta} f({U_{\varepsilon}(\tau)})\right\|_{L^{2}}+\left\|\mathcal{D}_{+}^{\beta}{P_\varepsilon}(\cdot,\tau)\right\|_{L^{2}}\Big)d\tau\\\nonumber
   &&\hspace{5mm}\\
   \nonumber
   &&\leq \|\frac{d}{dt^{2}}(S_{\alpha}^{})_{\varepsilon}(t)\mathcal{D}_{+}^{\beta}Q_{\varepsilon}\|_{L^{2}}\\\nonumber
  &&\hspace{5mm}+\frac{C}{\Gamma(\alpha-1)}\int_{0}^{t}(t-\tau)^{\alpha-2}
  \|f'(U_{\varepsilon}(\tau))\|_{L^{\infty}}\|\mathcal{D}_{+}^{\beta}\partial_{\tau}U_{\varepsilon}(\tau)\|_{L^{2}}d\tau\\
   \nonumber
   &&\hspace{5mm} + \frac{C}{\Gamma(\alpha-1)}\int_{0}^{t}(t-\tau)^{\alpha-2}
  \|\mathcal{D}_{+}^{\beta}f'(U_{\varepsilon}(\tau))\|_{L^{2}}\|\partial_{\tau}U_{\varepsilon}(\tau)\|_{L^{\infty}}d\tau\\
   \nonumber
  &&\hspace{5mm}+\frac{1}{\Gamma(\alpha-1)}\int_{0}^{t}(t-\tau)^{\alpha-2}
  \|\mathcal{D}_{+}^{\beta}\partial_{\tau}{P_\varepsilon}(\cdot,\tau)\|_{L^{2}}d\tau\\
  \nonumber
  &&\hspace{5mm} +\Big(C\|f'(Q_{\varepsilon})\|_{L^{\infty}}\|\mathcal{D}_{+}^{\beta}Q_{\varepsilon}\|_{L^{2}}
  +\|\mathcal{D}_{+}^{\beta}{P_\varepsilon}(\cdot,0))\|_{L^{2}}\Big)\frac{t^{\alpha-2}}{\Gamma(\alpha-1)}\\
  \nonumber
  &&\hspace{5mm}
  +C\int_{0}^{t}(t-\tau)^{2\alpha-3}\|\widetilde{A}_{\varepsilon}\|\,
  E_{\alpha,2\alpha-2}((t-\tau)^{\alpha}\|\widetilde{A}_{\varepsilon}\|)
  \,\|f'(U_{\varepsilon}(\tau))\|_{L^{\infty}}\|\mathcal{D}_{+}^{\beta}U_{\varepsilon}(\tau)\|_{L^{2}}d\tau\\
    \nonumber
  &&\hspace{5mm}+
  \int_{0}^{t}(t-\tau)^{2\alpha-3}\|\widetilde{A}_{\varepsilon}\|\,
  E_{\alpha,2\alpha-2}((t-\tau)^{\alpha}\|\widetilde{A}_{\varepsilon}\|)
  \,\|\mathcal{D}_{+}^{\beta}{P_\varepsilon}(\cdot,\tau)\|_{L^{2}}d\tau.
 \end{eqnarray}
Similarly to (\ref{34b}), based on assumptions for function $f'$  and
using the estimate (\ref{2alpha}), one gets
 \begin{equation*}
  \lim_{t\rightarrow0^{+}}\left\|\frac{\mathcal{D}_{+}^{\beta}\frac{d}{dt^{2}}U_{\varepsilon}(t,\cdot)}{t^{\alpha-2}}\right\|_{L^{2}}=C<+\infty,
\end{equation*}
 Thus, it follows that the property (\ref{26a}) is satisfied.

 Next, we prove the moderate bound for $\|^{C}\mathcal{D}_{t}^{\gamma}U_{\varepsilon}(t,\cdot)\|_{H^{\beta}},$ where $\gamma \in \{0,1,\alpha\},$ and also for $\|\frac{d}{dt^{2}}U_{\varepsilon}(t,\cdot)\|_{H^{\beta}}$. First, we prove the moderate bound for $\|^{C}\mathcal{D}_{t}^{\gamma}U_{\varepsilon}(t,\cdot)\|_{H^{\beta}},$
 considering the cases:
\begin{enumerate}
  \item $ Case \; \gamma=0$\\
  From the representation (\ref{7b}) and alternative integral representation given in Proposition \ref{intrepre1}, one has
  \begin{eqnarray}\label{13a}
  \nonumber
    \|U_{\varepsilon}(t)\|_{L^{2}}&\leq& \|(S_{\alpha})_{\varepsilon}(t)Q_{\varepsilon}\|_{L^{2}}\\\nonumber
&&+\int\limits_{0}^{t}(t-\tau)^{\alpha-1}\|
E_{\alpha,\alpha}((t-\tau)^{\alpha}\widetilde{A}_{\varepsilon})\|\cdot\|f(U_{\varepsilon})\|_{L^{2}}d\tau\\\nonumber
&&+\int\limits_{0}^{t}(t-\tau)^{\alpha-1}\|
E_{\alpha,\alpha}((t-\tau)^{\alpha}\widetilde{A}_{\varepsilon})\|\cdot\|{P_\varepsilon}(\cdot,\tau)\|_{L^{2}}d\tau\\\nonumber
  &\leq& \|(S_{\alpha})_{\varepsilon}(t)Q_{\varepsilon}\|_{L^{2}}\\\nonumber
  &&+\int\limits_{0}^{t}\frac{(t-\tau)^{\alpha-1}}{\alpha-1}
  E_{\alpha,\alpha-1}((t-\tau)^{\alpha}\|\widetilde{A}_{\varepsilon}\|)\|f(U_{\varepsilon})\|_{L^{2}}d\tau \\\nonumber
  &&+\int\limits_{0}^{t}\frac{(t-\tau)^{\alpha-1}}{\alpha-1}
  E_{\alpha,\alpha-1}((t-\tau)^{\alpha}\|\widetilde{A}_{\varepsilon}\|)\|{P_\varepsilon}(\cdot,\tau)\|_{L^{2}}d\tau .
  \end{eqnarray}

Using the estimate for $E_{\alpha,\alpha-1}$ given by (\ref{19}) one gets
\begin{eqnarray*}
  E_{\alpha,\alpha-1}(t^{\alpha}\|\widetilde{A}_{\varepsilon}\|) &\leq& C_{\alpha}(1+\|\widetilde{A}_{\varepsilon}\|^{(2-\alpha)/\alpha})
  (1+t^{2-\alpha})exp(t\|\widetilde{A}_{\varepsilon}\|^{1/\alpha}),
\end{eqnarray*}
so, if we introduce a notation
\begin{equation}\label{31c}
\widetilde{M}_{T}:=\sup_{t\in[0,T)}E_{\alpha,\alpha-1}(t^{\alpha}\|\widetilde{A}_{\varepsilon}\|),
\end{equation}
 then based on the well known properties of Landau symbols $o$ and $\mathcal{O},$ one gets the estimate
\begin{eqnarray}\label{32c}
\nonumber
  \widetilde{M}_{T} &\leq& C_{\alpha}\Big(1+o\Big(({\rm log}({\rm log}1/\varepsilon)^{\alpha-1})^{2-\alpha}\Big)\Big)\Big(1+T^{2-\alpha}\Big) \\\nonumber
   && \cdot {\rm exp}\Big(T\cdot o\Big({\rm log}({\rm log}1/\varepsilon)^{\alpha-1}\Big)\Big)\\\nonumber
   &\leq&  C_{\alpha}\Big(1+o\Big(({\rm log}{\rm log}1/\varepsilon)^{2-\alpha}\Big)\Big)\Big(1+T^{2-\alpha}\Big)\cdot {\rm exp}\Big(o\Big({\rm log}({\rm log}1/\varepsilon)^{\alpha-1}\Big)\Big)\\\nonumber
   &\leq& C_{\alpha}\Big(1+o\Big(({\rm log}1/\varepsilon)^{2-\alpha}\Big)\Big)\Big(1+T^{2-\alpha}\Big)\cdot o\Big(({\rm log}1/\varepsilon)^{\alpha-1}\Big)\\\nonumber
   &=& C_{\alpha}\Big(1+T^{2-\alpha}\Big)\Big(o\Big(({\rm log}1/\varepsilon)^{\alpha-1}\Big)+o({\rm log}1/\varepsilon)\Big)\\\nonumber
   &=& C_{\alpha}\Big(1+T^{2-\alpha}\Big)o({\rm log}1/\varepsilon)\\
   &=&\mathcal{O}({\rm log}1/\varepsilon).
\end{eqnarray}

Next, using (\ref{31c}) one gets the estimate
\begin{eqnarray*}
      \|U_{\varepsilon}(t)\|_{L^{2}} &\leq& \|(S_{\alpha})_{\varepsilon}(t)Q_{\varepsilon}\|_{L^{2}}
      +\frac{\widetilde{M}_{T}}{\alpha-1}\int\limits_{0}^{t}(t-\tau)^{\alpha-1}(\|f(U_{\varepsilon})\|_{L^{2}}+\|{P_\varepsilon}(\cdot,\tau)\|_{L^{2}})d\tau \\
       &\leq& \|(S_{\alpha})_{\varepsilon}(t)Q_{\varepsilon}\|_{L^{2}}+\frac{\widetilde{M}_{T}}{\alpha-1}\int\limits_{0}^{t}(t-\tau)^{\alpha-1}(C\|U_{\varepsilon}(\tau)\|_{L^{2}}+\|{P_\varepsilon}(\cdot,\tau)\|_{L^{2}})d\tau,
    \end{eqnarray*}
    and taking into account the relation (\ref{32c}),  the Gronwall's inequality gives the moderate bound for $\|U_{\varepsilon}(t)\|_{L^{2}}$.

    After differentiation of (\ref{7b}) with respect to left Liouville fractional derivative $\mathcal{D}_{+}^{\beta},$ using integral representation similarly to the one given in Proposition \ref{intrepre1}, we have $$\mathcal{D}_{+}^{\beta}U_{\varepsilon}(t)=(S_{\alpha})_{\varepsilon}(t)\mathcal{D}_{+}^{\beta}Q_{\varepsilon}
    +\int\limits_{0}^{t}(t-\tau)^{\alpha-1}E_{\alpha,\alpha}((t-\tau)^{\alpha}
    \widetilde{A}_{\varepsilon})(\mathcal{D}_{+}^{\beta}f(U_{\varepsilon})+\mathcal{D}_{+}^{\beta}{P_\varepsilon}(\cdot,\tau))d\tau,$$
    and similarly to the previous estimates one gets
    \begin{eqnarray*}
      \|\mathcal{D}_{+}^{\beta}U_{\varepsilon}(t)\|_{L^{2}} &\leq& \|(S_{\alpha})_{\varepsilon}(t)\mathcal{D}_{+}^{\beta}Q_{\varepsilon}\|_{L^{2}}\\
      &&+
    \int\limits_{0}^{t}\frac{(t-\tau)^{\alpha-1}}{\alpha-1}E_{\alpha,\alpha-1}((t-\tau)^{\alpha}\|\widetilde{A}_{\varepsilon}\|)
    \|\mathcal{D}_{+}^{\beta}(f(U_{\varepsilon}))\|_{L^{2}}d\tau\\
    &&+
    \int\limits_{0}^{t}\frac{(t-\tau)^{\alpha-1}}{\alpha-1}E_{\alpha,\alpha-1}((t-\tau)^{\alpha}\|\widetilde{A}_{\varepsilon}\|)
    \|\mathcal{D}_{+}^{\beta}{P_\varepsilon}(\cdot,\tau)\|_{L^{2}}d\tau\\
    &\leq& \|(S_{\alpha})_{\varepsilon}(t)\mathcal{D}_{+}^{\beta}Q_{\varepsilon}\|_{L^{2}}\\
    &&+\frac{\widetilde{M}_{T}}{\alpha-1}\int\limits_{0}^{t}(t-\tau)^{\alpha-1}\|f'(U_{\varepsilon}(\tau))\|_{L^{\infty}}\|\mathcal{D}_{+}^{\beta}U_{\varepsilon}(\tau)\|_{L^{2}}d\tau\\
    &&+\frac{\widetilde{M}_{T}}{\alpha-1}\int\limits_{0}^{t}(t-\tau)^{\alpha-1}\|\mathcal{D}_{+}^{\beta}{P_\varepsilon}(\cdot,\tau)\|_{L^{2}}d\tau.
    \end{eqnarray*}
    Taking into account assumptions for function $f',$ the moderate bound for  $\|\mathcal{D}_{+}^{\beta}U_{\varepsilon}(t)\|_{L^{2}}$ follows after applying the Gronwall's inequality.
   \item $Case \; \gamma=1$\\
   Using integral representation for (\ref{7b}), given by (\ref{35}), after differentiation  with respect to $t,$ one gets
   \begin{equation}\label{prvi}
     \frac{d}{dt}U_{\varepsilon}(t)=\frac{d}{dt}(S_{\alpha})_{\varepsilon}(t)Q_{\varepsilon}
    +\int\limits_{0}^{t}(t-\tau)^{\alpha-2}E_{\alpha,\alpha-1}((t-\tau)^{\alpha}
    \widetilde{A}_{\varepsilon})(f(U_{\varepsilon})+{P_\varepsilon}(\cdot,\tau))d\tau,
   \end{equation}
   and then the estimate
    \begin{eqnarray*}
      \|\frac{d}{dt}U_{\varepsilon}(t)\|_{L^{2}} &\leq& \|\frac{d}{dt}(S_{\alpha})_{\varepsilon}(t)Q_{\varepsilon}\|_{L^{2}} + \widetilde{M}_{T}\int\limits_{0}^{t}(t-\tau)^{\alpha-2}(C\|U_{\varepsilon}(\tau)\|_{L^{2}}+\|{P_\varepsilon}(\cdot,\tau)\|_{L^{2}})d\tau.
    \end{eqnarray*}
    Thus, the moderate bound for  $\|\frac{d}{dt}U_{\varepsilon}(t)\|_{L^{2}}$ immediately follows.

    Similarly, after differentiation of (\ref{prvi}) with respect to left Liouville fractional derivative $\mathcal{D}_{+}^{\beta},$ we have
    \begin{eqnarray}\label{drugi}
     \nonumber
      \|\mathcal{D}_{+}^{\beta}\frac{d}{dt}U_{\varepsilon}(t)\|_{L^{2}} &\leq& \|\frac{d}{dt}(S_{\alpha})_{\varepsilon}(t)\mathcal{D}_{+}^{\beta}Q_{\varepsilon}\|_{L^{2}} \\\nonumber
       && + \widetilde{M}_{T}\int\limits_{0}^{t}(t-\tau)^{\alpha-2}\|f'(U_{\varepsilon}(\tau))\|_{L^{\infty}}\|\mathcal{D}_{+}^{\beta}U_{\varepsilon}(\tau)\|_{L^{2}}d\tau\\
       && + \widetilde{M}_{T}\int\limits_{0}^{t}(t-\tau)^{\alpha-2}\|\mathcal{D}_{+}^{\beta}{P_\varepsilon}(\cdot,\tau)\|_{L^{2}}d\tau,
    \end{eqnarray}
    and again the moderate bound for $\|\mathcal{D}_{+}^{\beta}\frac{d}{dt}U_{\varepsilon}(t)\|_{L^{2}}$ immediately follows.
  \item $Case \; \gamma=\alpha$\\
 From (\ref{6b}) one gets $$\|^{C}\mathcal{D}_{t}^{\alpha}U_{\varepsilon}(t)\|_{L^{2}}\leq
  \|\widetilde{A}_{\varepsilon}U_{\varepsilon}(t)\|_{L^{2}}+\|f(U_{\varepsilon})\|_{L^{2}}+\|{P_\varepsilon}(\cdot,t)\|_{L^{2}},$$ and since $f$ is globally Lipschitz function and holds $f(0)=0,$
  one gets the moderate bound for $\|^{C}\mathcal{D}_{t}^{\alpha}U_{\varepsilon}(t)\|_{L^{2}}.$

 After differentiation of (\ref{6b}) with respect to left Liouville fractional derivative $\mathcal{D}_{+}^{\beta},$ we have
  \begin{eqnarray*}
    \|\mathcal{D}_{+}^{\beta}\;^{C}\mathcal{D}_{t}^{\alpha}U_{\varepsilon}(t)\|_{L^{2}} &\leq& \|\mathcal{D}_{+}^{\beta}(\widetilde{A}_{\varepsilon}U_{\varepsilon}(t))\|_{L^{2}}+\|\mathcal{D}_{+}^{\beta}(f(U_{\varepsilon}))\|_{L^{2}} +\|\mathcal{D}_{+}^{\beta}{P_\varepsilon}(\cdot,t))\|_{L^{2}}\\
     &\leq& C(\log1/\varepsilon)^{\alpha} \|U_{\varepsilon}(t)\|_{H^{\beta}}+\|f'(U_{\varepsilon}(t))\|_{L^{\infty}}\|\mathcal{D}_{+}^{\beta}U_{\varepsilon}(t)\|_{L^{2}}   \\
     &&+ \|\mathcal{D}_{+}^{\beta}{P_\varepsilon}(\cdot,t)\|_{L^{2}},
  \end{eqnarray*}
  and using similarly arguments one gets the moderate bound for $\|\mathcal{D}_{+}^{\beta}\;^{C}\mathcal{D}_{t}^{\alpha}U_{\varepsilon}(t)\|_{L^{2}}.$
\end{enumerate}

It remains to prove the moderate bound for $\|\frac{d}{dt^{2}}U_{\varepsilon}(t,\cdot)\|_{H^{\beta}}.$ This time, the moderate bound  is obtained using inequalities  (\ref{34b}) and (\ref{32c}), as well as previously shown moderate bounds for $\|U_{\varepsilon}(t,\cdot)\|_{H^{\beta}}$  and $\|\frac{d}{dt}U_{\varepsilon}(t,\cdot)\|_{H^{\beta}}.$

To prove that this solution is unique in Colombeau space $\mathcal{G}_{\alpha}([0,\infty):H^{\beta}(\mathbb{R})),$ suppose that there exist two solutions $U$ and $V$ to (\ref{6b}) and set $\xi_{\varepsilon}=U_{\varepsilon}-V_{\varepsilon}$. This difference satisfies
\begin{eqnarray}\label{11b}
 \nonumber
  ^{C}\mathcal{D}_{t}^{\alpha}\xi_{\varepsilon}(t) &=& \widetilde{A}_{\varepsilon}\xi_{\varepsilon}(t)+f(U_{\varepsilon})-f(V_{\varepsilon}) + \widetilde N_{\varepsilon}(t),\\
  \xi_{\varepsilon}(0) &=& \xi_{0\varepsilon}, \;\;(\xi_{t})_{\varepsilon}(0)=0,
\end{eqnarray}
where $\widetilde N_{\varepsilon}(t)\in \mathcal{N}_{\alpha}([0,\infty):H^{\beta}(\mathbb{R}))$ and $\xi_{0\varepsilon}\in \mathcal{N}_{\alpha}(H^{\beta}(\mathbb{R}))$. Then, solution of the problem (\ref{11b}) is given by
\begin{eqnarray}\label{7c}
\nonumber
  \xi_{\varepsilon}(t) &=& (S_{\alpha})_{\varepsilon}(t)\xi_{0\varepsilon}
+\int_{0}^{t}\,_{\tau}J_{t}(S_{\alpha})_{\varepsilon}(t-\tau)^{RL}\mathcal{D}_{\tau}^{2-\alpha}
(f(U_{\varepsilon})-f(V_{\varepsilon}))d\tau \\
   &&+\int\limits_{0}^{t}\,_{\tau}J_{t}(S_{\alpha})_{\varepsilon}(t-\tau)^{RL}\mathcal{D}_{\tau}^{2-\alpha}\widetilde N_{\varepsilon}(\tau)d\tau.
   \end{eqnarray}
Using alternative integral representation (\ref{7c}) becomes
\begin{eqnarray}\label{12b}
\nonumber
  \xi_{\varepsilon}(t) &=& (S_{\alpha})_{\varepsilon}(t)\xi_{0\varepsilon}
+\int\limits_{0}^{t}(t-\tau)^{\alpha-1}E_{\alpha,\alpha}((t-\tau)^{\alpha}\widetilde{A}_{\varepsilon})(f(U_{\varepsilon})-f(V_{\varepsilon}))d\tau \\
   &&+\int\limits_{0}^{t}(t-\tau)^{\alpha-1}E_{\alpha,\alpha}((t-\tau)^{\alpha}\widetilde{A}_{\varepsilon})\widetilde N_{\varepsilon}(\tau)d\tau,
   \end{eqnarray}
and further one gets following estimate
\begin{eqnarray*}
  \|\xi_{\varepsilon}(t)\|_{L^{2}} &\leq& \|(S_{\alpha})_{\varepsilon}(t)\xi_{0\varepsilon}\|_{L^{2}}\\ &&+\int\limits_{0}^{t}(t-\tau)^{\alpha-1}\|E_{\alpha,\alpha}((t-\tau)^{\alpha}\widetilde{A}_{\varepsilon})
  \|\cdot\|f(U_{\varepsilon})-f(V_{\varepsilon})\|_{L^{2}}d\tau\\
   && +\int\limits_{0}^{t}(t-\tau)^{\alpha-1}\|E_{\alpha,\alpha}((t-\tau)^{\alpha}\widetilde{A}_{\varepsilon})\|\cdot\|\widetilde N_{\varepsilon}(\tau)\|_{L^{2}}d\tau. \end{eqnarray*}
   Now, using $\widetilde{M}_{T}$ given by (\ref{31c}) and estimated by (\ref{32c}), as well as the fact that $f$ is Lipschitz function with respect to $u$, after applying the  Gronwall's inequality one gets the $\mathcal{N}$-bound for $\|\xi_{\varepsilon}(t)\|_{L^{2}}$.

   From (\ref{11b}) it follows the estimate $$\|^{C}\mathcal{D}_{t}^{\alpha}\xi_{\varepsilon}(t)\|_{L^{2}}\leq
   \|\widetilde{A}_{\varepsilon}\xi_{\varepsilon}(t)\|_{L^{2}}+\|f(U_{\varepsilon})-f(V_{\varepsilon}) \|_{L^{2}}+\|\widetilde N_{\varepsilon}(t)\|_{L^{2}},$$
   and then the $\mathcal{N}$-bound for $\|^{C}\mathcal{D}_{t}^{\alpha}\xi_{\varepsilon}(t)\|_{L^{2}}$ is easy to obtain.

Next,
\begin{eqnarray*}
  \|\mathcal{D}_{+}^{\beta}\xi_{\varepsilon}\|_{L^{2}} & \leq &
  \|(S_{\alpha})_{\varepsilon}(t)\mathcal{D}_{+}^{\beta}\xi_{0\varepsilon}\|_{L^{2}} \\
  & + & \int\limits_{0}^{t}\|(S_{\alpha})_{\varepsilon}(t-s)
  \mathcal{D}_{+}^{\beta}(f(U_{\varepsilon}(\cdot,s))-f(V_{\varepsilon}(\cdot,s)))\|_{L^{2}}ds \\ & + & \int\limits_{0}^{t}\|(S_{\alpha})_{\varepsilon}(t-s)\mathcal{D}_{+}^{\beta}N_{\varepsilon}(s)\|_{L^{2}}.
\end{eqnarray*}
From the mean value theorem we have
$$
f(U_{\varepsilon}(\cdot,s))-f(V_{\varepsilon}(\cdot,s))=(U_{\varepsilon}(\cdot,s)-V_{\varepsilon}(\cdot,s))f'(\theta
  U_{\varepsilon}(\cdot,s)+(1-\theta)V_{\varepsilon}(\cdot,s)),\; 0< \theta <1.
$$
Introduce the notation
\begin{equation*}
F(\cdot,s)=f'(\theta U_{\varepsilon}(\cdot,s)+(1-\theta)V_{\varepsilon}(\cdot,s)).
\end{equation*}
From the approximate Leibniz rule for the Liouville fractional derivative (Proposition \ref{23}) we get
\begin{eqnarray*}
  && \|\mathcal{D}_{+}^{\beta}(f(U_{\varepsilon}(\cdot,s))-f(V_{\varepsilon}(\cdot,s)))\|_{L^{2}}=
  \|\mathcal{D}_{+}^{\beta}[F(\cdot,s)(U_{\varepsilon}(\cdot,s)-V_{\varepsilon}(\cdot,s))]\|_{L^{2}} \\
   && \hspace{2mm} \leq C\|F\|_{L^{\infty}}\|\mathcal{D}_{+}^{\beta}(U_{\varepsilon}(s)-V_{\varepsilon}(s))\|_{L^{2}}
  +C\|\mathcal{D}_{+}^{\beta}F\|_{L^{2}}\|U_{\varepsilon}(s)-V_{\varepsilon}(s)\|_{L^{\infty}}\\
   && \hspace{2mm} \leq C\|F\|_{L^{\infty}}\|\mathcal{D}_{+}^{\beta}\xi_{\varepsilon}(s)\|_{L^{2}}
   +C\|\mathcal{D}_{+}^{\beta}F\|_{L^{2}}\|\xi_{\varepsilon}(s)\|_{H^{\beta}}.
\end{eqnarray*}
 so the $\mathcal{N}$-bound for $\|\mathcal{D}_{+}^{\beta}\xi_{\varepsilon}(t)\|_{L^{2}}$ follows after applying the Gronwall's inequality.

$\mathcal{N}$-bound for $\|\frac{d}{dt}\xi_{\varepsilon}(t)\|_{L^{2}}$ is possible to obtain after differentiation of (\ref{12b}) with respect to $t.$  In that away one gets
\begin{eqnarray*}
  &&\|\frac{d}{dt}\xi_{\varepsilon}(t)\|_{L^{2}} \leq \|\frac{d}{dt}(S_{\alpha})_{\varepsilon}(t)\xi_{0\varepsilon}\|_{L^{2}}\\ &&\hspace{10mm}+\int\limits_{0}^{t}(t-\tau)^{\alpha-2}\|E_{\alpha,\alpha-1}((t-\tau)^{\alpha}\widetilde{A}_{\varepsilon})
  \|\cdot\|f(U_{\varepsilon})-f(V_{\varepsilon})\|_{L^{2}}d\tau\\
   && \hspace{10mm}+\int\limits_{0}^{t}(t-\tau)^{\alpha-2}\|E_{\alpha,\alpha-1}((t-\tau)^{\alpha}\widetilde{A}_{\varepsilon})\|\cdot\|\widetilde N_{\varepsilon}(\tau)\|_{L^{2}}d\tau,
\end{eqnarray*}
and similarly to (\ref{prvi}) it follows the $\mathcal{N}$-bound. Also, after differentiation of (\ref{12b}) first with respect to $t,$ and then with respect to left Liouville fractional derivative $\mathcal{D}_{+}^{\beta},$  similarly to (\ref{drugi}), one gets the $\mathcal{N}$-bound for $\|\mathcal{D}_{+}^{\beta}\frac{d}{dt}\xi_{\varepsilon}(t)\|_{L^{2}}.$

 Differentiation of (\ref{11b})  with respect to left Liouville fractional derivative $\mathcal{D}_{+}^{\beta},$ implies
 \begin{eqnarray*}
   \|\mathcal{D}_{+}^{\beta}\;^{C}\mathcal{D}_{t}^{\alpha}\xi_{\varepsilon}(t)\|_{L^{2}} &\leq& \|\mathcal{D}_{+}^{\beta}(\widetilde{A}_{\varepsilon}\xi_{\varepsilon}(t))\|_{L^{2}}+\|\mathcal{D}_{+}^{\beta}(f(U_{\varepsilon})-f(V_{\varepsilon})) \|_{L^{2}} \\
    && +\|\mathcal{D}_{+}^{\beta}\widetilde N_{\varepsilon}(t)\|_{L^{2}},
 \end{eqnarray*}
 and $\mathcal{N}$-bound for $\|\mathcal{D}_{+}^{\beta}\; ^{C}\mathcal{D}_{t}^{\alpha}\xi_{\varepsilon}(t)\|_{L^{2}}$ is obtained using similar arguments as in the previous cases.

  $\mathcal{N}$-bound for $\|\frac{d}{dt^{2}}\xi_{\varepsilon}(t)\|_{H^{\beta}}$ is possible to obtain firstly differentiating (\ref{7c}) twice with respect to $t$, and then differentiating the newly equation with respect to left Liouville fractional derivative. In both cases, the corresponding $\mathcal{N}$-bound is obtained using integral representation (\ref{35c}), and applying the procedure similarly to the one using in (\ref{34b}) and (\ref{35b}).

 From all previous estimates it follows $\xi_{\varepsilon}:=U_{\varepsilon}-V_{\varepsilon}\in \mathcal{N}_{\alpha}([0,\infty):H^{\beta}(\mathbb{R})),$ i.e. the solution of fractional Cauchy problem  (\ref{6b}), represented by (\ref{7b}), is unique.
 \end{proof}

\begin{remark}\label{remark}
If $\widetilde{A}\in \mathcal{SG}(H^{2}(\mathbb{R}))$ is an operator of $h_{\varepsilon}$-type with $h_{\varepsilon}= o\Big((\log(\log1/\varepsilon)^{\alpha-1})^{\alpha}\Big)$, similarly one can prove that solution to (\ref{6b}) is also represented by (\ref{7b}) and this unique solution belongs to $\mathcal{G}_{\alpha}([0,\infty):H^{2}(\mathbb{R})).$
\end{remark}

\begin{definition}
The solution $U$ of the problem (\ref{6b}) introduced in Theorem \ref{14} is called generalized solution of the equation
\begin{equation*}
  ^{C}\mathcal{D}_{t}^{\alpha}U(t)=\widetilde{A}U(t)+f(U)+P(\cdot,t)
\end{equation*}
with generalized operators.
\end{definition}

\begin{remark}
 The generalized operator $\widetilde{A}$ satisfying properties (i) and (ii) can be obtained, for instance, by regularization of space fractional or integer order derivatives appearing in the operator $A.$
 For details we refer to \cite{Japundzic}.
 %The role of properties (i) and (ii) for various types of fractional equations will be seen in the next section.
\end{remark}

\section{Application to solving time and time-space fractional wave equations}

As said in the Introduction, in case of equations wave equations (\ref{fwe-nonfrac}) and (\ref{fwe}) one can obtain the approximate operator $\tilde A$ for a given (integer or fractional) differential operator $A$ by regularizing the derivative. In examples below, one can prove that the operators $A$ and $\widetilde A$ are $L^2$-associated (for details we refer \cite{Japundzic}).

\subsection{\textbf{Time fractional wave equation}}

 Let $1<\alpha<2$ and instead of the Cauchy problem (\ref{fraceq}) with integer space derivative, i.e. \begin{eqnarray*}
 \nonumber
  ^{C}\mathcal{D}_{t}^{\alpha}u(x,t) &=& \lambda(x)
\partial_{x}^{2}u(x,t)+f(u(x,t))+P(x,t), \hspace{2mm} t>0, \hspace{2mm} x\in \mathbb{R},\\
  u(0) &=& u_{0}, \;\;u_{t}(0) = 0,
\end{eqnarray*}
let us consider the corresponding approximate problem:
 \begin{eqnarray*}
 \nonumber
  ^{C}\mathcal{D}_{t}^{\alpha}U_{\varepsilon}(x,t) &=& \widetilde{A}U_{\varepsilon}(x,t)+f(U_{\varepsilon}(x,t))+{P_\varepsilon}(x,t), \hspace{2mm} t>0, \hspace{2mm} x\in \mathbb{R},\\
  Q_{\varepsilon} &=& u_{0 \varepsilon}, \;\;(u_{t})_{\varepsilon}(0) = 0,
\end{eqnarray*}
where the operator $\widetilde{A}\in \mathcal{SG}(H^{2}(\mathbb{R}))$ is represented by the nets of operators
\begin{eqnarray*}
 & &
\widetilde{A}_{\varepsilon}:H^{2}(\mathbb{R})\rightarrow H^{2}(\mathbb{R}), \\ & &
\widetilde{A}_{\varepsilon}U_{\varepsilon}=\lambda_{\varepsilon}(x)(
\partial_{x}^{2}U_{\varepsilon}*\phi_{h_{\varepsilon}}),
\end{eqnarray*}
such that
$\lambda_{\varepsilon}\in
H^{2}(\mathbb{R}),$
$\|\lambda_{\varepsilon}\|_{H^{2}(\mathbb{R})}=\mathcal{O}\Big((\log(\log1/\varepsilon)^{\alpha-1})^{\alpha/2}\Big),$
  $\phi_{h_{\varepsilon}}(x)=h_{\varepsilon}\phi(xh_{\varepsilon}),$ where $h_{\varepsilon}= o\Big((\log(\log1/\varepsilon)^{\alpha-1})^{\alpha/5}\Big),$ $\phi\in C_{0}^{\infty}(\mathbb{R}),$
  $\phi(x)\geq 0 $ and $\int\phi(x)dx=1.$

Then, from Theorem \ref{14} and Remark \ref{remark} it follows that the unique solution is represented by (\ref{7b}) and belongs to Colombeau space $\mathcal{G}_{\alpha}([0,\infty):H^{\beta}(\mathbb{R})), \beta \in (1,2).$

\subsection{ \textbf{Time-space fractional wave equation}}

Let us consider the approximate Cauchy problem for the time-space fractional wave equation, i.e.
\begin{equation*}
  ^{C}\mathcal{D}_{t}^{\alpha}u(t)=\widetilde{A}_{\beta}u(t)+f(u)+P(\cdot,t),
\end{equation*}
where  $1<\alpha<2,$ $1<\beta<2,$ the operator $\widetilde{A}_{\beta}\in \mathcal{SG}(H^{2}(\mathbb{R}))$ is represented by the nets of operators
\begin{eqnarray*}
 & &
(\widetilde{A}_{\beta})_{\varepsilon}:H^{2}(\mathbb{R})\rightarrow H^{2}(\mathbb{R}), \\ & &
(\widetilde{A}_{\beta})_{\varepsilon}U_{\varepsilon}=\lambda_{\varepsilon}(x)(
\mathcal{D}_{+}^{\beta}U_{\varepsilon}*\phi_{h_{\varepsilon}}),
\end{eqnarray*}
$\mathcal{D}_{+}^{\beta}$ is the left Liouville fractional derivative of order $\beta$  on the whole axis $\mathbb{R}$ given by
$$(\mathcal{D}_{+}^{\beta}u)(x)=\frac{1}{\Gamma(2-\beta)}\left(\frac{d}{dx}\right)^{2}
\int\limits_{-\infty}^{x}\frac{u(\xi)}{(x-\xi)^{\beta-1}}d\xi,$$
$\lambda_{\varepsilon}\in
H^{2}(\mathbb{R})$ and  $\phi_{h_{\varepsilon}}(x)$ satisfy the same properties as in the case of time fractional wave equation.

Again, from Theorem \ref{14} and Remark \ref{remark} it follows that the unique solution is represented by (\ref{7b}) and belongs to Colombeau space $\mathcal{G}_{\alpha}([0,\infty):H^{\beta}(\mathbb{R})).$ The same result holds if instead of left $\beta$th Liouville fractional derivative in the fractional operator $\widetilde{A}_{\beta},$ $1<\beta<2,$ one uses right $\beta$th Liouville fractional derivative or Riesz $\beta$th fractional derivative (for details we refer \cite{Japundzic}).

\remark The procedure proposed in the paper can also be applied to solve the problem (\ref{fwe}) with inhomogeneous second initial condition, i.e.
\begin{eqnarray*}
\nonumber
  ^{C}{\mathcal{D}}_{t}^{\alpha}u(x,t) &=& \lambda(x){\mathcal{D}}_{x}^{\beta}u(x,t)+f(u(x,t)), \,\, (x,t)\in \mathbb{R}\times \mathbb{R}_{+}, \\
  u(x,0)&=& u_{0}(x), \;\;\; u_{t}(x,0)= u_{1}(x).
\end{eqnarray*}
In that case the representation formula for approximate solution would be given by
\begin{equation*}
  u(t)=S_{\alpha}(t)u_{0}+JS_{\alpha}(t)u_{1}
+\int_{0}^{t}\;_{\tau}J_{t}S_{\alpha}(t-\tau)^{RL}\mathcal{D}_{\tau}^{2-\alpha}f(u(\tau))d\tau,
\end{equation*}
where  $S_{\alpha}$  is a Colombeau uniformly continuous solution operator generated by $\widetilde{A}$ which approximates originate differential operator.
The assertion concerning existence and uniqueness of approximate solution, as well as applications to fractional wave equations with variable coefficients continue to be valid.

\section*{Acknowledgements}
This work was supported by
Science Fund of the Republic of Serbia, GRANT No TF C1389-YF (“FluidVarVisc”)


\begin{thebibliography}{9}

\bibitem{Adams}R. Adams,  \textit{Sobolev Spaces}, Academic Press, New York, 1975.

\bibitem{Atanackovic} T.M. Atanackovi\'{c}, S. Pilipovi\'{c}, B. Stankovi\'{c}, D. Zorica, \textit{Fractional Calculus with Applications in Mechanics, Vibrations and Diffusion Processes}, ISTE, London, John Wiley and Sons, New York, (2014).

\bibitem{Atanackovic1} T.M. Atanackovi\'{c}, S. Pilipovi\'{c}, B. Stankovi\'{c}, D. Zorica, \textit{Fractional calculus with Application in Mechanics: Wave Propagation, Impact and Variational Principles}, ISTE, London, John Wiley and Sons, New York, (2014).

\bibitem{Bazhlekova} E. Bazhlekova, \textit{Fractional Evolution Equations in Banach Spaces}, PhD thesis, Eindhoven University of Technology, (2001).

\bibitem{Biagoni} H.A. Biagioni, \textit{A Nonlinear Theory of Generalized Functions}, Lect. Notes Math. 1421, Springer Verlag, Berlin, (1990).

\bibitem{Colombeau} J. F. Colombeau, \textit{Elementary Introduction to New Generalized Functions}, North Holland, Amsterdam, (1985).

\bibitem{Japundzic} M. Japund\v{z}i\'{c}, D. Rajter-\'{C}iri\'{c}, \textit{Reaction-advection-diffusion equations with space fractional derivatives and variable coefficients on infinite domain}, Fract. Calc. Appl. Anal., \textbf{18} (4) (2015), 911--950.

 \bibitem{EJDE} M. Japund\v{z}i\'{c}, D. Rajter-\'{C}iri\'{c}, \textit{Generalized uniformly continuous solution operators and inhomogeneous fractional evolution equations with variable coefficients}, Electron. J. Differ. Eq., \textbf{293} (2017), 1--24.

 \bibitem{Applicable} M. Japund\v{z}i\'{c}, D. Rajter-\'{C}iri\'{c}, \textit{Approximate solutions of time and time-space fractional wave equations with variable coefficients}, Appl. Anal., \textbf{97} (9) (2018), 1565--1590.

 \bibitem{Jennings} G. I. Jennings, \textit{Efficient Numerical Methods for Water Wave Propagation in Unbounded Domains}, PhD thesis, University of Michigan, (2012).

\bibitem{Kilbas} A. Kilbas, H. Srivastava and J. Trujillo, \textit{Theory and Applications of Fractional Differential Equations}, North-Holland, Amsterdam, (2006).

\bibitem{Nedeljkovb} M. Nedeljkov, S. Pilipovi\'c and D. Scarpal$\acute{e}$zos, \textit{The Linear Theory of Colombeau Generalized Functions}, Pitman Res. Not. Math., Longman, Essex, (1998).

\bibitem{Rajter} M. Nedeljkov, D. Rajter-\'Ciri\'c, Semigroups and PDEs with perturbations. In: Delcroix, A., Hasler, M.,Marti, J.-A., Valmorin, V. (eds.) Nonlinear Algebraic Analysis and Applications. Proc. of ICGF 2000, (2004), 219--227.

\bibitem{Nedeljkov} M. Nedeljkov, D. Rajter-\'Ciri\'c, \textit{Generalized uniformly continuous semigroups and semilinear
    hyperbolic systems with regularized
    derivatives}, Monatsh. Math., \textbf{160} (2010), 81--93.

\bibitem{Oberguggenberger} M. Oberguggenberger, \textit{Multiplication of Distributions and Applications to Partial
    Differential Equations}, Pitman Res. Not. Math. Vol.
    259, Longman Sci. Techn., Essex, (1992).

\bibitem{Pazy} A. Pazy, \textit{Semigroups of Linear Operators and Applications to Partial Differential Equations},
    Springer-Verlag, New York, (1983).

\bibitem{Podlubny} I. Podlubny, \textit{Fractional Differential Equations}, Academic Press, San Diego, (1999).

\bibitem{Samko} S. Samko, A. Kilbas, O. Marichev, \textit{Fractional Integrals and Derivatives}, Gordon and Breach Science Publishers, Amsterdam, (1993).

\bibitem{Umarov} S. Umarov, E. Saydamatov, \textit{A fractional analog of the Duhamel principle}, Fract. Calc. Appl. Anal., \textbf{9} (1) (2006), 57--70.

\bibitem{Webb} M. Webb,  \textit{Analysis and Approximation of a Fractional Differential
Equation}, Master’s Thesis, University of Oxford, Hilary Term, 2012.

\end{thebibliography}
\end{document}